\documentclass[10pt]{article}

\usepackage{amsbsy,amsfonts,amsmath}
\usepackage{caption,graphicx, color}
\usepackage[text={6.75in,9.5in},centering,letterpaper]{geometry}
\usepackage[comma,square,sort&compress,numbers]{natbib}
\usepackage{paralist}

\usepackage[hypertexnames=false,colorlinks=true,linkcolor=blue,citecolor=blue]{hyperref}

\newtheorem{Lemma}{Lemma}[section]

\newtheorem{Definition}[Lemma]{Definition}
\newtheorem{Hypothesis}[Lemma]{Hypothesis}

\newtheorem{Theorem}{Theorem}

\newenvironment{Proof}[1][.]%
 {\begin{trivlist}\item[]\textbf{Proof#1 }}%
 {\hspace*{\fill}$\rule{0.3\baselineskip}{0.35\baselineskip}$\end{trivlist}}

 {\begin{trivlist}\item[]\textbf{Acknowledgments }}{\end{trivlist}}

\makeatletter\@addtoreset{equation}{section}\makeatother

\def\Re{\mathop\mathrm{Re}\nolimits}    

\newcommand{\C}{\mathbb{C}}             
\newcommand{\N}{\mathbb{N}}             
\newcommand{\R}{\mathbb{R}}             
\newcommand{\Z}{\mathbb{Z}}             

\newcommand{\rmD}{\mathrm{D}}
\newcommand{\rmi}{\mathrm{i}}           
\newcommand{\rmW}{\mathrm{W}}

\newcommand{\prwt}[1]{P_{\mathrm{wt}}^\mathrm{#1}}

\newcommand{\bigo}{\mathcal{O}}
\newcommand{\expsmall}{\mathcal{O}(\mathrm{e}^{-\eta L})} 
\newcommand{\Rg}{\mathrm{Rg}}

\newcommand{\bigzero}{\mbox{\normalfont\Large\bfseries 0}}

\newcommand{\changed}[1]{#1}

\begin{document}

\title{Homoclinic snaking of contact defects in
reaction-diffusion equations}

\author{%
Timothy Roberts\footnote{Corresponding Author: email timothyr@uchicago.edu} \\
Department of Statistics\\
University of Chicago\\
Chicago, IL~60637, USA
\and%
Bj\"orn Sandstede\\
Division of Applied Mathematics\\
Brown University\\
Providence, RI~02912, USA
}

\date{\today}
\maketitle

\begin{abstract}
We apply spatial dynamical-systems techniques to prove that certain spatiotemporal patterns in reversible reaction-diffusion equations undergo snaking bifurcations. That is, in a narrow region of parameter space, countably many branches of patterned states coexist that connect at towers of saddle-node bifurcations. Our patterns of interest are contact defects, which are 1-dimensional time-periodic patterns with a spatially oscillating core region that at large distances from the origin in space resemble pure temporally oscillatory states and arise as natural analogues of spiral and target waves in one spatial dimension. We show that these solutions lie on snaking branches that have a more complex structure than has been seen in other contexts. In particular, we predict the existence of families of asymmetric traveling defect solutions with arbitrary background phase offsets, in addition to symmetric standing target and spiral patterns. We prove the presence of these additional patterns by reconciling results in classic ODE studies with results from the spatial-dynamics study of patterns in PDEs and using geometrical information contained in the stable and unstable manifolds of the background wave trains and their natural equivariance structure.
\end{abstract}


\section{Introduction} 

Pattern formation is ubiquitous in myriad applications from physics, chemistry, and biology. A principal example of these are spiral and target waves which have garnered much attention in both mathematical and applied literature in recent years \cite{glass1996, ertl1991, lee1996, cross1993}. In this paper we will present new bifurcation results for one-dimensional spiral and target patterns in reaction-diffusion equations via snaking bifurcations, \changed{which were obtained in the thesis \cite{roberts2024} of the first author}.

\begin{figure}
    \centering
    \includegraphics[height=6cm]{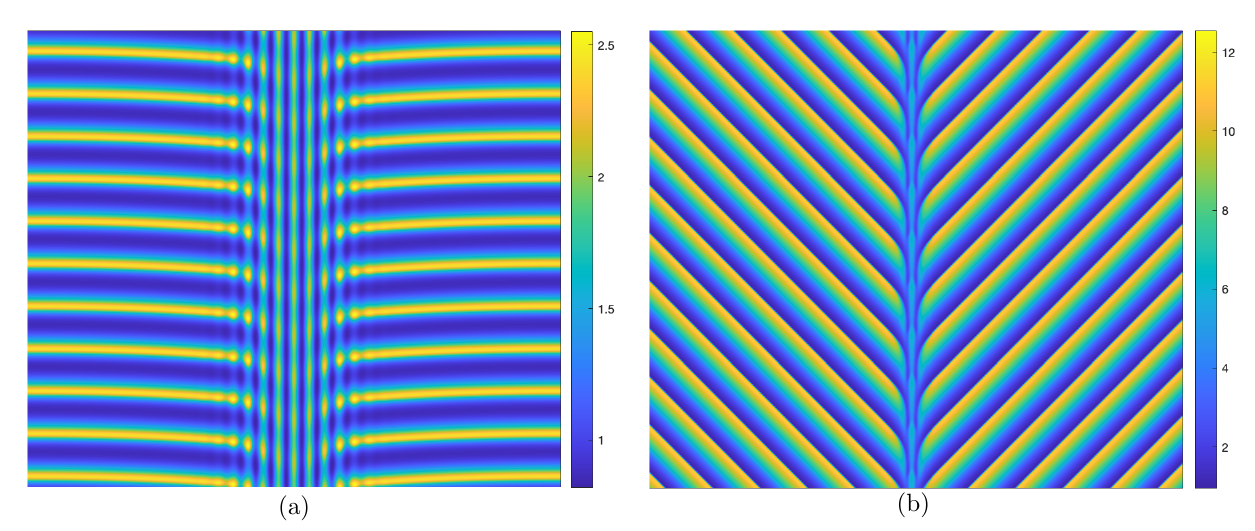}
    \caption{Space-time plots of example defects solutions for the Brusselator  \eqref{eqn:bruss_1} for different choices of parameter values. The horizontal axis denotes space and the vertical axis depicts time. (a) Contact defect. (b) Source defect.}
    \label{fig:defect}
\end{figure}

Our patterns of interest are defect solutions \cite{sandstede2004}, \changed{which are time-periodic in the laboratory or an appropriate comoving coordinate frame and which consist of a spatial core region that mediates between two traveling waves with periodic profiles (which we refer to as wave trains) in the far field as illustrated in Figure~\ref{fig:defect}.} Such solutions have been observed experimentally (see, for instance, \cite{perraud1993}) and studied both numerically and analytically in many applications. In particular, they provide a one-dimensional analogue of spatiotemporal spiral and target waves in reaction-diffusion systems. In this context, spiral and target waves are differentiated by the phase offset between their background oscillations, where spirals have a phase difference of half a wavelength, while targets have in-phase oscillations.

\changed{Defect solutions in reaction-diffusion systems can be classified into four distinct types \cite{sandstede2004}. Here, we review two of them, namely source and contact defects, which are relevant for this paper. Source defects are locally unique and converge exponentially in space to the asymptotic wave trains in the far field. They are active defects so that localized perturbations are transported away from the core towards the far field, and their temporal frequency, computed in a comoving coordinate frame if the defect is traveling, is selected nonlinearly via a Rankine--Hugoniot condition \cite[(1.8)]{sandstede2004}. Figure~\ref{fig:defect}(b) contains a space-time plot of a stationary source defect.}

\changed{Of primary interest to us are \textit{contact defects}. Stationary contact defects are asymptotic to the same spatially homogeneous, temporal oscillation in the far field as visible in Figure~\ref{fig:defect}(a). Their temporal frequency is identical to the frequency of the asymptotic oscillation, and there is therefore no nonlinear frequency selection by the contact defect itself. The slopes of the curves of constant phase in the space-time plot decay like $1/|x|$ (see \eqref{e:velocity} below for an explanation) instead of exponentially as for source defects. Stationary contact defects are generically accompanied by traveling contact defects that approach the same small-wavenumber wave train as $x\to\pm\infty$ and travel with the group velocity of the asymptotic wave train \cite{sandstede2004,sandstede2004a}.}

\begin{figure}
    \centering
    \includegraphics[width=0.9\linewidth]{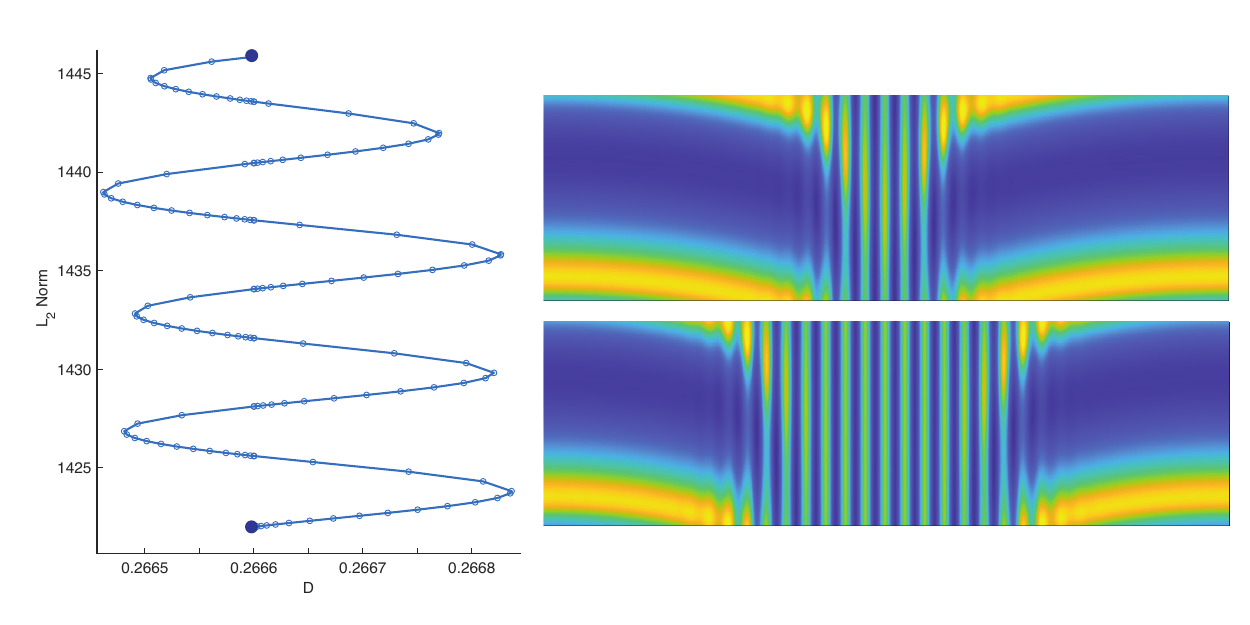}
    \caption{A reproduction of the numerical observations of \cite{tzou2013}. Branches of contact defects are computed on finite domains and continued in system parameters. The defects are seen to snake since stripes on the temporally homogeneous core are gained as one traverses the curve. Example space-time plots for defects are shown for a single temporal period.}
    \label{fig:tzou_numerics}
\end{figure}

Numerical studies by Tzou et al.~\cite{tzou2013} provided evidence for interesting new bifurcation processes for defect solutions in reaction-diffusion systems, and especially the Brusselator. They numerically observed that symmetric contact defects exist in the neighborhood of a codimension-2 Turing--Hopf bifurcation and that they appear to snake as the diffusion ratio is swept through a narrow region of parameter space. The contact defects themselves are a mixed state where the core region is a \changed{stationary spatially-periodic} Turing pattern and \changed{the far-field patterns at either side are given by a spatially homogeneous, temporally oscillatory state arising from the Hopf mode}. Their studies indicate that in a narrow region of parameter space there exists a tower of saddle-node bifurcations joining a countably infinite number of branches of contact defects. Moreover, as the tower is traversed, the inner core region is allowed to grow. A reproduction of their numerics is shown in Figure~\ref{fig:tzou_numerics}.

Similar behavior has been studied extensively in the context of exponentially localized stationary solutions for PDEs. Work in this direction began through studying stationary solutions in applications to mechanics, fluids, and optics \cite{champneys1998}. Homoclinic and heteroclinic connections were studied in the unfolding of a codimension-2 reversible Hamiltonian--Hopf (or reversible 1-1 resonance) bifurcation \cite{woods1999}. Here it was found that patterns of interest resulted from complicated intersections of stable and unstable manifolds called heteroclinic tangles. However, the full picture of these bifurcation diagrams was not understood until the work of Burke and Knobloch \cite{burke2007, burke2007a} who identified not only the snakes of symmetric branches but also the ladder rungs of asymmetric branches joining them --- hence the name snakes and ladders bifurcation. Work by Beck et al.\ \cite{beck2009} established a rigorous theory from which one can predict the qualitative features of the snaking diagram, and move beyond the realm of conservative and low dimensional systems. 

\begin{figure}
    \centering
    \includegraphics[width=0.9\textwidth]{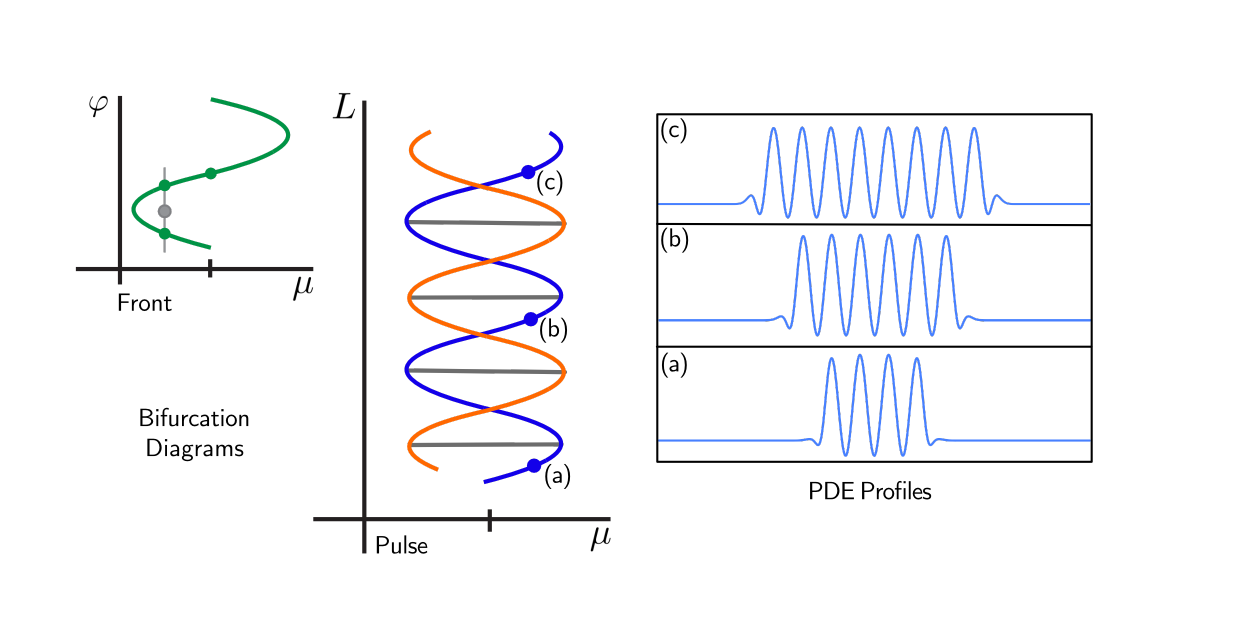}
    \caption{Construction of the snaking diagram for the Swift-Hohenberg equation. We assume that the bifurcation diagram of the front is periodic in a natural variable $\varphi$. The bifurcation diagram of exponentially localized pulse solutions is then formed by two periodic extensions of the front diagram offset by half a period. Asymmetric branches are formed at the midpoint between distinct fronts at the same values of the bifurcation parameter $\mu$. As you move up the diagram, the core region width $L$ increases.}
    \label{fig:homoclinic_snaking}
\end{figure}

The theory of homoclinic snaking in the case of localized patterns can be understood through the lens of homoclinic and heteroclinic bifurcations. \changed{One assumes the presence of fronts: a stationary (time-independent) solution that is asymptotic in backward space to the rest state zero and in forward space to a spatially periodic pattern such as a Turing pattern. By reversibility this implies the existence of a commensurate back solution found by reversing the spatial coordinate that is asymptotic in backward space to the periodic state and in forward space to zero.} One then asks whether, if we were to cut these solutions at the right places, they could be glued together to form an exponentially localized solution with a spatially periodic core region. In the conservative setting, for example the quadratic-cubic or cubic-quintic Swift--Hohenberg equation \cite{beck2009}, the answer is affirmative. Moreover, if you assume that the fronts exist and continue in a bifurcation parameter (and make a global assumption about the shape of the corresponding front bifurcation diagram --- that it is an isola), then you can qualitatively predict the resulting bifurcation diagram. There are two branches of symmetric solutions each of which winds back and forth through the bifurcation interval at a countable succession of saddle nodes, as pictured in Figure~\ref{fig:homoclinic_snaking}. Moving up through each saddle node bifurcation the inner core region grows. These symmetric branches are connected at periodic intervals by branches of asymmetric solutions like rungs of a ladder, and the shape of the snaking diagram can be predicted using only the bifurcation diagram of the front solutions (up to exponentially small errors). This prediction procedure is illustrated in Figure~\ref{fig:homoclinic_snaking}.

Here we present a rigorous framework for understanding the existence and continuation of contact defects in general reaction-diffusion systems with reversible symmetry. We show that these patterns undergo snaking bifurcations and the snaking bifurcation diagrams found have oppositely winding symmetric branches connected by asymmetric branches, analogous to the localized case. However, upon adding the asymmetric solutions, we find that bifurcation diagrams are markedly different and more complicated than previous results suggest. Instead of finding a discrete set of asymmetric solutions for each parameter value, we find 1-parameter families parameterized by a phase offset in the background oscillations of the defect solutions so that the asymmetric branches interpolate between the symmetric spiral and target waves already found. We describe this in more detail in the following section. \changed{We emphasize that we cannot prove that the Brussellator system (or any other pattern-forming reaction-diffusion system) satisfies the hypotheses we need for our construction. However, our assumptions consist of nondegeneracy conditions, which we expect to hold for typical systems.}

\subsection{Main Results}

While our results, as stated in Section~\ref{section:intuition_infinitedim}, are valid for a large family of reaction-diffusion equations, we can give a concrete example in terms of the Brusselator as our motivating equation:
\begin{align}
    \begin{split}
        U_t &= DU_{xx} + E - (B+1)U + VU^2,  \\
        V_t &= V_{xx} + BU - VU^2.
    \end{split} \label{eqn:bruss_1}
\end{align}
In this instance we can state our results concretely as follows. First we make assumptions on the geometry of the phase space: the precise formulation of the genericity conditions we outline below is involved, and a full account of these is given in Section~\ref{section:set_up_and_hyp}. \changed{Throughout, we set $S^1:=\R/2\pi\Z$.}

\begin{Hypothesis} \label{hyp:Brusselator}
Suppose there exists a smooth function $(z,\kappa^*):S^1 \rightarrow J\times K$ defined on intervals $J \subseteq \R$ and $K \subseteq \R^+$ with non-empty interior such that for every $\mu=z(s)$ in $J$, there exists a front solution to the Brusselator for diffusion ratio $D = \mu$ which is asymptotic in forward space to the Turing pattern with wavenumber $\kappa^*(s)$ and asymptotic in backward space to the spatially homogeneous oscillation with temporal frequency $\omega_0(s)\neq0$. Furthermore, we assume that this bifurcation family of fronts satisfies genericity conditions corresponding to the transversality and general position of the stable and unstable manifolds of the Turing patterns and the homogeneous oscillations, respectively.
\end{Hypothesis}

Assuming these conditions hold we prove the following two theorems.

\begin{Theorem}[Snaking of Symmetric Defects] \label{thm:bruss_symmetric_snaking}
\changed{Assume Hypothesis~\ref{hyp:Brusselator} holds, then there are constants $N\in\N$ and $\eta>0$ so that the following holds. For each $n>N$, each $\varphi_0\in\{0,\pi\}$, and each $s\in S^1$, there exists a unique symmetric target contact defect and a unique symmetric spiral contact defect solution at the parameter value $\mu=z(s)+\bigo(\mathrm{e}^{-\eta n})$ with core width $L=(2\pi n+\varphi_0-s)/\kappa^*(s)+\bigo(\mathrm{e}^{-\eta n})$ and temporal frequency $\omega=\omega_0(s)+\bigo(\mathrm{e}^{-\eta n})$, where the Turing patterns in the core of the defect have wavenumber $\kappa=\kappa^*(s)+\bigo(\mathrm{e}^{-\eta n})$ and a minimum at $x=0$ for $\varphi_0=0$ and a maximum for $\varphi_0=\pi$.} Here, uniqueness is up to space and time translations.
\end{Theorem}

\changed{In particular, symmetric target and spiral contact defects each arise as separate one-parameter families where the parameter $\mu$ and the Turing plateau width $L$ vary as functions of the branch parameter $s$. Next, we consider asymmetric contact defects. We denote by $\theta$ the asymptotic phase difference between the asymptotic wave train profiles at $x=\pm\infty$, so that $\theta=0$ for target and $\theta=\pi$ for spiral contact defects. We remark that it was shown in \cite[\S6.2]{sandstede2004} that this phase difference is well defined (see also \S\ref{s:geom}).}

\begin{Theorem}[Snaking of Asymmetric Defects] \label{thm:bruss_asymmetric_snaking}
\changed{Assume Hypothesis~\ref{hyp:Brusselator} holds, then there are constants $N\in\N$ and $\eta > 0$ so that the following holds. Choose $n>N$, $\theta\in S^1$, and a pair $(s_1,s_2)\in S^1\times S^1$ that satisfies $z(s_1)=z(s_2)$ and $z^\prime(s_1)z^\prime(s_2)\neq0$. Then there exists a unique contact defect (up to translations in space and time) at the parameter value $\mu=z(s_1)+\bigo(\mathrm{e}^{-\eta n})$ with asymptotic phase difference $\theta$ and core width $L$ given by}
\begin{equation}
L = \frac{2\pi n-(s_1+s_2)}{\kappa^*(s_1)+\kappa^*(s_2)} + \bigo \left(|\kappa^*(s_1)-\kappa^*(s_2)|+\mathrm{e}^{-\eta n} \right). \label{eqn:asym_core_width}
\end{equation}
\changed{These defects move with wave speed $c_\mathrm{d}=\bigo(\frac{1}{n}|\kappa^*(s_1)-\kappa^*(s_2)|+\mathrm{e}^{-\eta n})$, have temporal frequency $\omega=\omega_0(s_1)+\bigo(c_\mathrm{d}^2)$, and connect background wave trains with spatial wavenumber $k=\bigo(c_\mathrm{d}^2)$.}
\end{Theorem}

\begin{figure}
    \centering
    \includegraphics[width=0.9\linewidth]{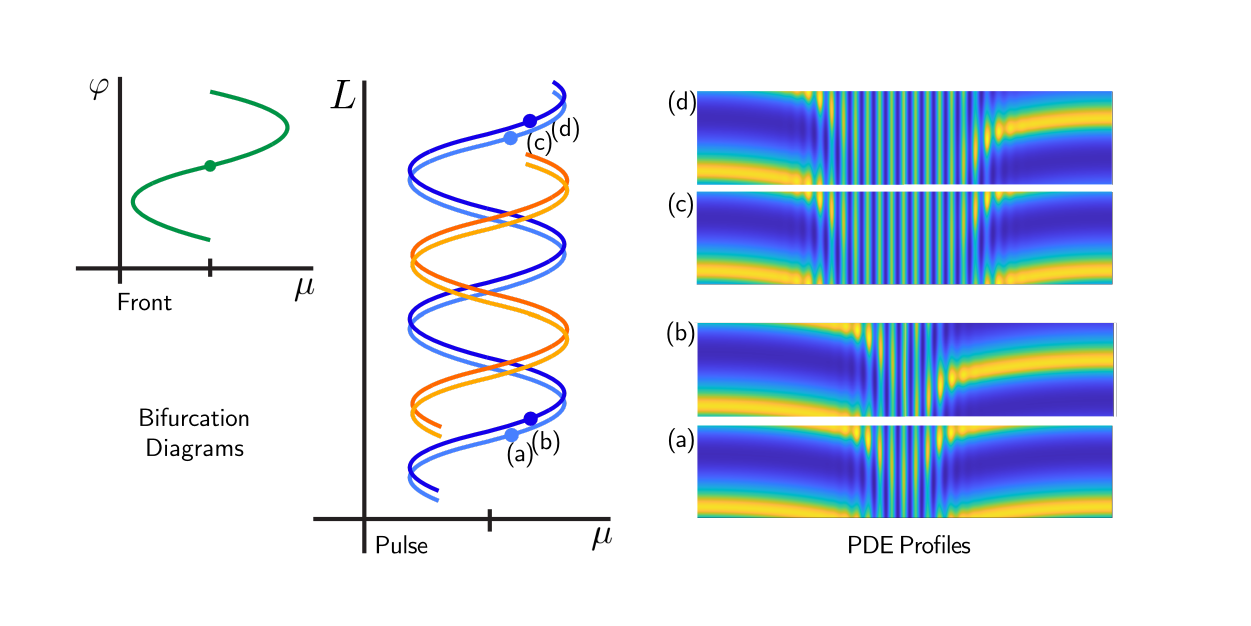}
    \caption{Qualitative description of the  snaking diagram for symmetric defects in the Brusselator. We assume the front connections have a periodic bifurcation diagram, displayed in green. To leading order, the symmetric bifurcation diagram is formed as a periodic tiling of the front bifurcation diagram. \changed{Two of the branches correspond to target patterns with, respectively, a minimum and a maximum at their center (see (a) and (c)), while the other two branches are spiral patterns with, respectively, a minimum and a maximum at their center (see (b) and (d))}.}
    \label{fig:prediction_sym}
\end{figure}

\begin{figure}
    \centering
    \includegraphics[width=0.75\linewidth]{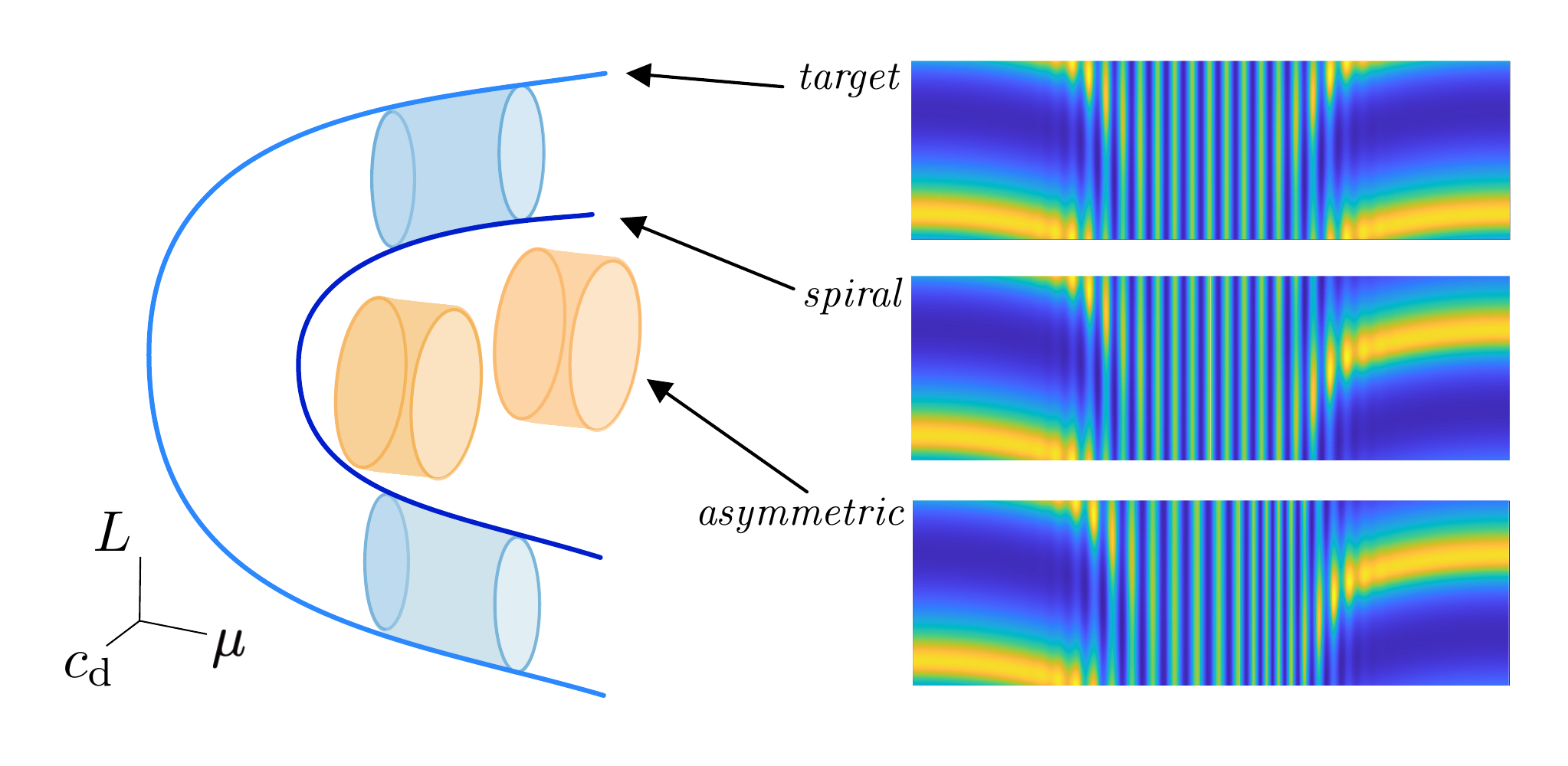}
    \caption{Qualitative description of the full snaking diagram in the Brusselator, where we plot an additional third dimension corresponding to the defect wave speed, $c_\mathrm{d}$. We again assume the front bifurcation diagram found in Figure~\ref{fig:prediction_sym} and zoom in on one saddle-node for clarity. \changed{In blue are the two symmetric defect bifurcation curves corresponding to targets and spirals respectively}. The blue tubular families have symmetric cores interpolating between the two stationary defect families \changed{by sweeping out all possible background phase offsets in between}. The orange families correspond to defects with asymmetric cores. The position of the asymmetric branches depends on the difference in selected Turing wavenumbers.}
    \label{fig:prediction_asym}
\end{figure}

These results correspond to the predicted bifurcation diagrams presented in Figures~\ref{fig:prediction_sym} and~\ref{fig:prediction_asym}. The results for symmetric defects are very similar to case of exponentially localized solutions: we find towers of saddle-nodes joining branches of symmetric defects whose shape traces out the periodic tiling of the front bifurcation diagram. The integer $n$ indexes the copy of the front diagram and $s$ the position along that diagram. As one climbs the tower, moving across saddle-node bifurcations, Turing stripes are added to the core region. There are four distinct symmetric snaking branches: two correspond to target patters (they have in-phase background oscillations), and two correspond to spiral waves (they have anti-phase background oscillations). Within these pairs, the patterns on the two oppositely winding branches differ in their symmetry at the origin, with one branch having a minimum at $x=0$ and the other a maximum.

\changed{In contrast, the asymmetric results are far more complex for defects than for exponentially localized stationary solutions. Indeed, there are now multiple distinct ways to construct asymmetric defects. Firstly, we can break the symmetry of the background oscillations, while keeping $s_1=s_2$ so that the core is fixed. For each defect solution on the symmetric branches, there is a 1-parameter family of defects whose cores are effectively the same but that have out-of-phase background oscillations}. These families come in circles parameterized by the phase offset $\theta$, which can take on all possible values in $[0,2\pi]$. This comes at the cost of both allowing the new asymmetric defects to move at a small but possibly nonzero wave speed $c_\mathrm{d}$ in the laboratory frame and changing the Turing wavenumber by a small amount. Moreover, Theorem~\ref{thm:asymmetric_contact} shows that these perturbations will be exponentially small in the length of the core Turing region, i.e.\ $c_\mathrm{d}=\bigo(\mathrm{e}^{-\eta n})$. Extending the bifurcation diagram to include a third dimension corresponding to the wave speed $c_\mathrm{d}$ (see Figure~\ref{fig:prediction_asym}), the bifurcation diagram now has tubular branches of solutions of which the symmetric target and spiral defects are two distinguished 1-parameter subfamilies which travel at wave speed zero. 

Alternatively, one can form asymmetric solutions by breaking the symmetry of the Turing-core, and possibly also the symmetry of the background oscillations. \changed{By choosing $s_1\neq s_2$ we find} defect solutions whose Turing spatial wavenumber is modulated. Moreover, these solutions also have out-of-phase background oscillations which produce more tubular families in the extended bifurcation diagrams seen in Figure~\ref{fig:prediction_asym}. In comparison to the former asymmetric solutions, these core-asymmetric defects will have larger wave speeds which decay algebraically as the core width grows with the rate determined by the difference $\Delta\kappa$ of the selected Turing wavenumbers so that $c_\mathrm{d}=\bigo(|\Delta\kappa|/n)$.

The existence of asymmetric solutions with arbitrary background phase offsets is a very surprising result in the context of spatial dynamical systems. Indeed to our knowledge, such defect solutions have not been found in any other context. This is because, in typical scenarios, there is a very weak (exponentially small, beyond-all-orders) coupling between the background oscillations which excludes the possibility of arbitrary phase offsets and only allows target and spiral defects \cite{chapman2009, kozyreff2006}. We find that this coupling is not present in the case of contact defects due to additional degrees of freedom introduced by additional center directions and algebraic convergence of the defects to their asymptotic oscillations. We discuss the intuition and underlying mechanisms behind the formation of these defects in Section~\ref{section:intuition_infinitedim}. However, it must be stressed that this is an essential feature for the snaking of contact defects and must be present in order to reconcile results from spatial dynamics with the theory of homoclinic snaking. We also see from this discussion that we should not think of these snaking contact defects as simply localized solutions on a non-uniform background. Indeed, it is only through a proper appreciation of their identity as defect solutions that we are able to state, prove and understand our results on these snaking bifurcations.

\subsection{Paper Outline}

We begin in Section~\ref{sec:main_results} by sketching the intuition behind homoclinic snaking in the case of localized rolls for scalar conservative systems. This is the original context for homoclinic snaking and will provide us with a foundation to understand the more complex situation of spatiotemporal patterns in reaction-diffusion systems. We will then discuss how this simpler picture needs to be updated in our new context. In particular, we discuss why it is essential that we find our additional families of asymmetric defect solutions by considering the new geometric picture provided by spatial dynamics. 

In Section~\ref{section:spatial_dynamics} we briefly summarize the background required to understand the full spatial dynamics picture, with a focus on the main players of our analysis: homogeneous oscillations, Turing patterns and defect solutions.

Section~\ref{section:set_up_and_hyp} starts the rigorous set up for our analysis including all required assumptions and a full discussion of our choices in describing the transversality and general position hypotheses required for the proofs. \changed{We end this section by restating Theorems~\ref{thm:bruss_symmetric_snaking} and~\ref{thm:bruss_asymmetric_snaking} in their full generality.} The proofs of \changed{the general} Theorems~\ref{thm:symmetric_contact} and~\ref{thm:asymmetric_contact} and detailed in Sections~\ref{section:proof_symmetric_snaking} and~\ref{section:asymmetric_snaking} respectively.

Finally we will discuss some outstanding questions about snaking of defect solutions in Section~\ref{section:discussion}. Here we will describe snaking for other types of defects and how that might change the associated bifurcation diagrams. We also discuss ongoing attempts to numerically approximate and continue contact defects in the Brusselator.

\section{Intuition and the Geometry of Phase Space} \label{sec:main_results}

\subsection{Finite Dimensional Analogue and Intuition}
\begin{figure}
    \centering
    \includegraphics[height=7cm]{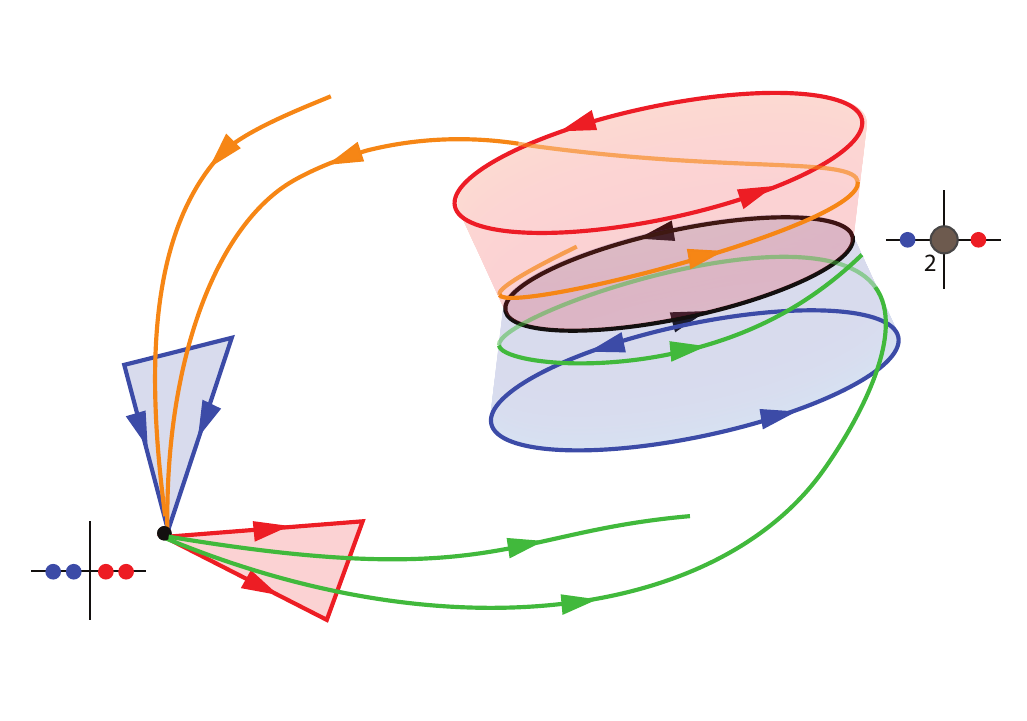}
    \caption{A schematic diagram of the fronts and backs in the scalar Swift--Hohenburg equation. In black are the rest state 0 and the periodic orbit $u_\mathrm{Tur}$ along with there stable and unstable manifolds in blue and red. The front and back are colored green and orange respectively.}
    \label{fig:schematic_SH}
\end{figure}

The term snaking bifurcation comes from studies of exponentially localized solutions in reversible ODEs. Since we borrow and extend this framework to the infinite dimensional case, we begin by describing the essential details in this simpler setting. As a model for the snaking process consider the scalar quadratic-cubic Swift--Hohenberg equation
\begin{align}
    U_t = -(1+\partial_x^2)^2U -\mu U + \nu U^2 -U^3,
\end{align}
although the exact details of the equation are less important than its properties. We wish to look for exponentially localized solutions which are stationary in time and asymptotic to the rest state zero. After setting $U_t = 0$ the task becomes locating homoclinic connections to the origin in a four dimensional phase space with dynamics
\begin{align}
    u_x = f(u; \mu). \label{eqn:SH_first_order}
\end{align}
where $u = (U,U_x,U_{xx},U_{xxx})$. We can immediately verify the following facts for a non-empty open interval of values $J$ for $\mu>0$. Equation \eqref{eqn:SH_first_order} is both conservative and reversible. Here reversible means that there exists a $\Z_2$-symmetry operation which, combined with a reversal of time, leaves the dynamics unchanged. The trivial rest state at zero is a saddle point with spectrum consisting of four real eigenvalues -- two stable and two unstable. There also exists a reversible periodic orbit, $\gamma(x;\mu)$, which has two Floquet exponents at zero (corresponding the trivial exponent and one associated with conservation) and two hyperbolic exponents one stable and one unstable. Reversibility of the periodic orbit requires the $\Z_2$-symmetry maps the periodic orbit to itself. A schematic of the geometry in phase space is sketched in Figure~\ref{fig:schematic_SH}.

The key property leading to the existence of localized patterns is the assumed existence of fronts. We assume that there exists a heteroclinic connection between the rest state zero in backward $x$ and the periodic orbit in forward $x$. While such a connection can be found numerically using approximation techniques, no formal proof of these structures currently exists, and any attempt at such a proof would require current advances in rigorous and validated numerics \cite{vandenberg2020,berg2023}. However, their presence is not unexpected and can be justified by dimension counting arguments. Looking inside the 3-dimensional zero energy level set, we have a 2-dimensional unstable manifold $\rmW^\mathrm{u}(0)$ and a 2-dimensional center-stable manifold $\rmW^\mathrm{cs}(\gamma(\cdot;\mu))$ which will have a 1-dimensional intersection generically, see Figure~\ref{fig:schematic_SH}. Meaning that if we have such an non-empty intersection for one value of $\mu$ this persists on a non-empty open interval.

By virtue of the reversibility, we also have a back solution. That is, a heteroclinic connection between the periodic orbit $\gamma$ in backward $x$ and the rest state in forward $x$. So we can naturally ask whether we can find a nearby homoclinic solution to zero that first follows the front into a neighborhood of the periodic before leaving to follow the back and returning to the origin. Again by dimension counting, we see that this constitutes an intersection between the two 2-dimensional manifolds $\rmW^\mathrm{u}(0)$ and $\rmW^\mathrm{s}(0)$ inside the 3-dimensional zero energy level set. Such an intersection is 1-dimensional generically, so it is reasonable to expect that such a homoclinic connection exists and can be continued in the bifurcation parameter $\mu$. 

Indeed, this turns out to be the case and, for fixed $\mu \in J$, you can find a discrete set of such homoclinic orbits indexed by an integer $n$ (for $n$ sufficiently large) corresponding to the number of revolutions the homoclinic connection makes while in transit of the periodic orbit. As already described these homoclinics are transverse and so persist under small variations in $\mu$. To determine the global behavior of the bifurcations, i.e. how these countably many branches connect together as we vary $\mu$, we need to make global assumptions on the nature of the fronts. In particular, we assume that the bifurcation diagram of the fronts is periodic (and satisfying non-degeneracy conditions), in which case we observe snaking. Moreover, the snaking diagram can be qualitatively and quantitatively predicted (up to exponentially small terms) using the bifurcation diagram of the heteroclinic fronts.

It is also possible to construct asymmetric solutions using the global assumptions on the heteroclinic fronts. According to this assumption, we can see that for each value of $\mu$ (away from the saddle-nodes) there are in fact two pairs of fronts and backs. Thus, we can also ask for homoclinic connections that follow one front and the other back (or vice versa). Doing so we find additional asymmetric homoclinic branches bifurcating at pitchfork bifurcations from the symmetric branches (which will be generic as they result from breaking the reversible symmetry of the system). Continuing these branches globally, you find that that they in fact connect the two oppositely winding symmetric branches like rungs of a ladder. Again there is a simple algebraic condition using the bifurcation diagram of the fronts that, qualitatively and quantitatively, describes the behavior of these branches up to exponentially small errors. 

\changed{Moving from the conservative, reversible Swift--Hohenberg equation to the non-conservative, reversible Brusselator system, we need to consider whether the lack of a conservative structure will impact our arguments.} By considering the geometry already set up, we can see that this is not the case. If we ignore the conservative structure of the equation, then we see that a generic intersection of the 2-dimensional unstable manifold of the origin and the 3-dimensional center-stable manifold of the periodic orbit will be 1-dimensional. So even without the conservative structure, the front and back are generic connections and will have the same existence and robustness properties that we have assumed already. The difference is that to make those connections we have to use the second center-direction of the periodic orbit, which in this case parameterizes the energy level sets nearby the periodic orbit. When not in the conservative setting, reversibility will again produce this additional direction and it can be used to make the same argument as above assuming one is able to account for the additional degree of freedom.

\subsection{Infinite Dimensional Problem} \label{section:intuition_infinitedim}

We can now extend this approach to the infinite dimensional case by treating the defects as heteroclinic connections between Turing and Hopf modes.

We consider reaction-diffusion equations of the form
\begin{align}
    u_t = D u_{xx} + f(u), \label{eqn:reaction-diffusion}
\end{align}
where $u \in \R^2$ is the dependent variable, $D$ is a diagonal matrix with strictly positive entries and $f$ is some suitable nonlinear function.

Our main tool is to think of  defect solutions as connecting orbits in a spatial dynamical system. In \eqref{eqn:reaction-diffusion}, we set $w = \partial_x u$ to find
\begin{align}
\begin{split}
    u_x &= w, \\ \label{eqn:spatial_dynamics}
    w_x &= D^{-1}(u_t - f(u)),
\end{split}
\end{align}
a dynamical system on the phase space $X = H^{1/2}_{\mathrm{per}}( S^1,\R^2)\times L^2_{\mathrm{per}}( S^1,\R^2)$ where $ S^1:=\R/2\pi\Z$. \changed{While this system is not well-posed as an initial value problem, we can construct} exponential dichotomies and invariant manifolds to the Turing and Hopf modes (see for instance \cite{sandstede2004})\footnote{\changed{The theory in \cite{sandstede2004} requires that the leading-order operator $\mathcal{A}=\begin{pmatrix} 0 & 1\\ \partial_t & 0\end{pmatrix}$ has bounded spectral projections. On the span of $\mathrm{e}^{\rmi t}$, $\mathcal{A}$ is represented by the matrix $\begin{pmatrix} 0 & 1\\ \rmi n & 0\end{pmatrix}$ and the coordinate change $\begin{pmatrix} \sqrt{n} & 0\\ 0 & 1\end{pmatrix}$, which reflects the $H^{1/2}_\mathrm{per}\times L^2_\mathrm{per}$ norm, transforms this matrix into $\sqrt{n}\begin{pmatrix} 0 & 1\\ \rmi & 0\end{pmatrix}$, which has bounded spectral projections. This would not work on spaces other than $H^{1/2}_\mathrm{per}\times L^2_\mathrm{per}$.}}. Taking this perspective, defect solutions become heteroclinic connections between equilibria which lie near heteroclinic fronts and backs.

\begin{figure}
    \centering
    \includegraphics[height=7cm]{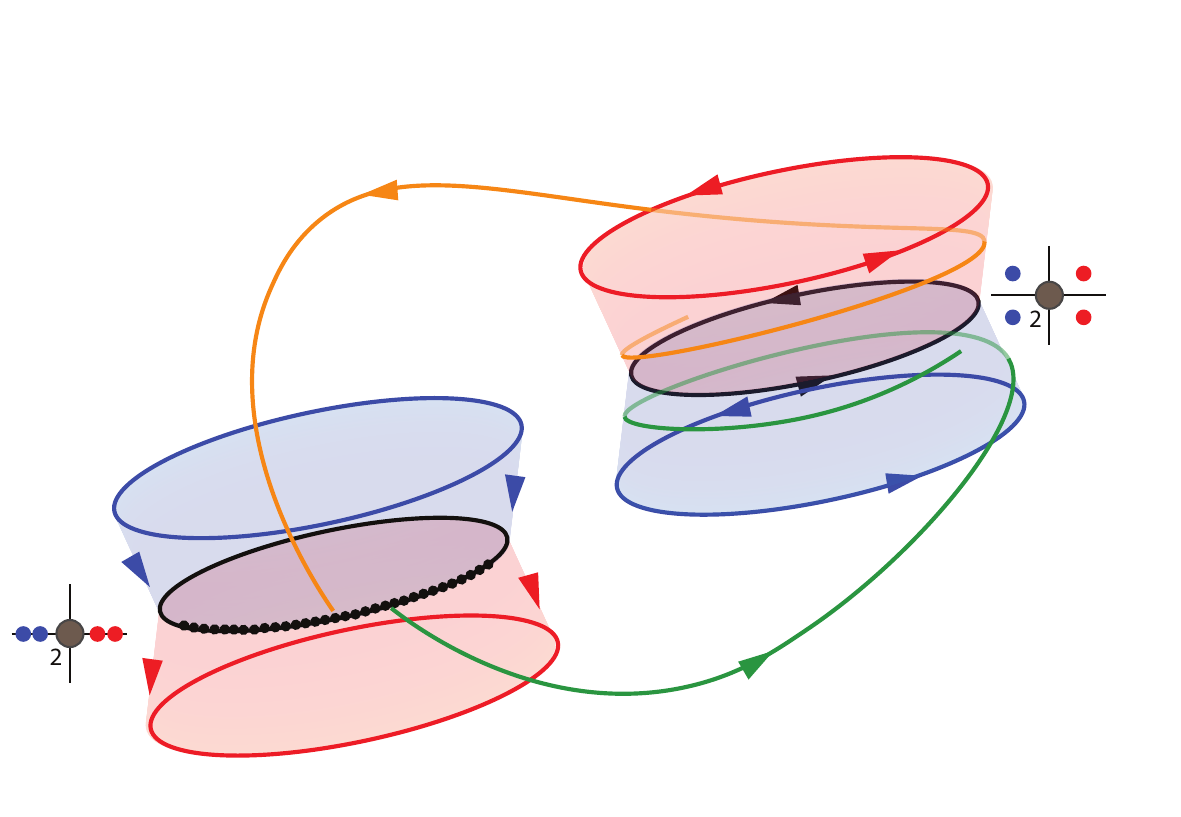}
    \caption{A schematic diagram of the fronts and backs in the reaction-diffusion case when looking for contact defects. In black are the rest state equilibria (Hopf modes) and a single periodic orbit (Turing pattern) $u_\mathrm{Tur}$ along with their stable and unstable manifolds in blue and red. Fronts are connections between the center-unstable manifold of the equilibria and the stable manifold of the periodic orbit; while backs are the reversed connections. In the symmetric case, fronts and backs select the same Turing pattern.}
    \label{fig:phase_space_contact}
\end{figure}

Using the spatial coordinate $x$ as the independent variable, we are able to describe (at lease schematically) the geometry of the problem in the same terms as the localized rolls. The asymptotic wave trains become circles of relative equilibria, topologically $ S^1$, with an $ S^1$-symmetry given by translations in time. Moreover from \cite{sandstede2004} we know precisely the spectrum of the background wave trains depending on their group velocities. For contact defects these circles have a two dimensional center subspace and all other spectrum bounded away from the imaginary axis. 

Simultaneously the Turing pattern becomes a periodic orbit sitting in the subspace of temporally constant functions. We also know that generically Turing patterns (viewed as periodic orbits of a reversible system) come in a 1-parameter family parameterized by their (spatial) \changed{wavenumber} \cite{vanderbauwhede1992, devaney1977}. This is reflected in their spectrum as a second Floquet exponent at zero. So we get an invariant cylinder of patterns representing candidate core regions for our defect solution.

We then assume that we have heteroclinic connections joining the background oscillations and the Turing patterns (and vise versa) which we will call fronts (and backs) from now on. At this point we are in the set up of homoclinic snaking. Defect solutions correspond to homoclinic connections in close proximity to the heteroclinic cycle of fronts and backs. However, the geometry of phase space is significantly more complicated than in the finite-dimensional case.

To build some intuition, we first look at a simpler situation where we imagine truncating the dimensions of our problem to keep only the bare essentials, see Figure~\ref{fig:phase_space_contact}. For our purposes this will mean considering the 6-dimensional subspace consisting of the six weakest directions. Intuitively speaking, assuming no higher codimension bifurcations occur, all important information is contained in these central dimensions. When it comes to actually proving our results, we must be careful as the problem is truly infinite-dimensional and we cannot hope to reduce the problem to a finite-dimensional system in general. However, the intuition of this model is essential in formulating the transversality and general position assumptions in our rigorous framework. In particular, now that we are in finite dimensions we are able to employ dimension counting arguments. 

First we observe that we expect the existence of a 1-parameter family of heteroclinic fronts. Indeed we have a 4-dimensional center-unstable manifold of the homogeneous oscillation intersecting the 4-dimensional center-stable manifold of the Turing patterns in $\R^6$. Such an intersection generically has dimension 2. This is essential as our PDE has translational symmetry in both time and space. Thus, we have a 1-parameter family of fronts where we expect the family to be parameterized by time translations, and that our front connections are generic and robust. At the same time, our previous dimension counting argument shows we should expect that defect solutions also come in 1-parameter families as both the center-stable and center-unstable manifolds of the homogeneous oscillations are 4-dimensional. 

\begin{figure}
    \centering
    \includegraphics[height=7cm]{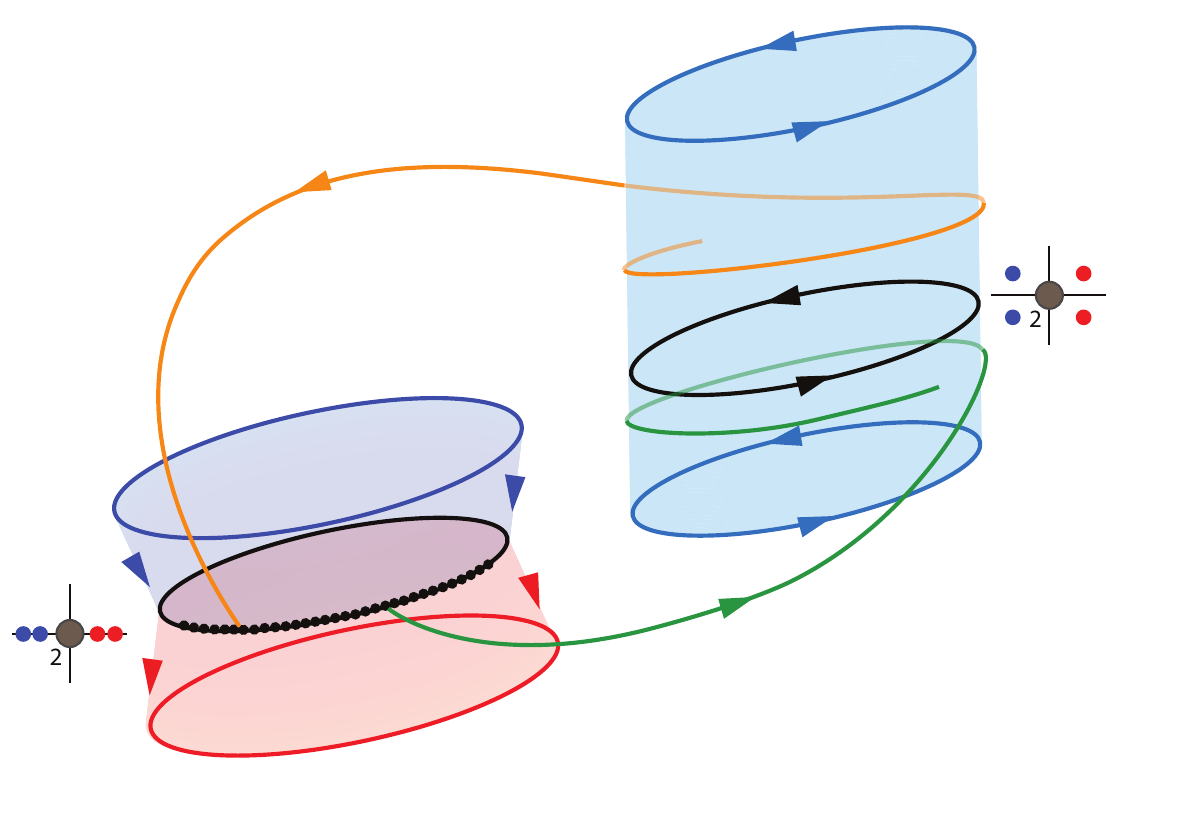}
    \caption{A schematic diagram of the fronts and backs in the reaction-diffusion case when looking for asymmetric contact defects. Here the fronts and backs select Turing patterns with different spatial wavenumbers. In black are the rest state equilibria (Hopf modes) and a single periodic orbit (Turing pattern) $u_\mathrm{Tur}$ along with their stable and unstable manifolds in blue and red. Representative fronts and backs are colored green and orange respectively.}
    \label{fig:phase_space_contact_asymm}
\end{figure}

We can then imagine trying to find asymmetric connections between two distinct fronts as in the Swift--Hohenberg case. However, here we find a problem. Generically fronts at different points on the bifurcation diagram will have core regions with different wavenumbers. As demonstrated above, our generic intersection for the fronts includes the center direction corresponding to the Turing wavenumber, meaning we should expect the fronts to have different core wavenumbers. For an illustration of this see Figure~\ref{fig:phase_space_contact_asymm}. The problem is that in the stationary laboratory frame, the dynamics on this cylinder are trivial, consisting of invariant circles of fixed wavenumber. So there is no way to move from a front with one wavenumber to a back with a different wavenumber. We have to move to a new comoving frame for some small, to be determined, wave speed $c_\mathrm{d}$ in order to make such a connection. This means that our generic intersection will have another additional dimension, not corresponding time translation. 

Instead this corresponds to a secondary operation which we will call shear translations. Suppose we take a symmetric defect and alter it by cutting it down the line $x=0$ and shifting only one half of the pattern forward in time by a phase \changed{$\theta$}. Looking at the resulting defect, the core region is essentially unchanged -- it is still a temporally homogeneous Turing pattern of the correct wavelength. It is natural to ask if this is again a solution. The answer is no and the resulting defect is no longer a solution, but it will be very close to one which moves at a very small {\it but nonzero} wave speed. These solutions correspond to the additional dimension in our family of solutions for contact defects. 

From a spatial dynamical systems viewpoint this is both unexpected and a critical feature. It is unexpected because in most applications one does not find defect solutions with arbitrary phase offsets. The typical scenario is that there is a weak coupling between the background oscillations that pin their phases to either be in-phase or out-of-phase by exactly half a wavelength, thereby producing only target and spiral defects \cite{kozyreff2006, chapman2009}. In this case, we actually find that there are solutions corresponding to each phase offset which interpolate between the spiral and target states.

On the other hand, this observation is also critical in making these results consistent with the robustness properties of contact defects. When $c$ becomes non-zero, we perturb the dynamics of the asymptotic homogeneous oscillations in addition to the family of Turing patterns. As a result, by moving into a comoving frame, we must now allow $\omega$ to vary so that we still have a contact defect. See Section~\ref{section:spatial_dynamics} for a discussion of this phenomenon. Consequently, the background states become wave trains with a well-defined but small spatial wavenumber $k$ which will be parameterized by the phase offset. This is important because \cite{sandstede2004} showed that contact defects are generically found in 1-parameter families parameterized by their far-field spatial wavenumbers, indicating that we should expect that there are many more solutions than were found in the numerical experiments by \cite{tzou2013} and that the background oscillations should not all be pure homogeneous oscillations. In their numerics, Tzou et al. compute only symmetric contact defects and did not find any asymmetric defects. From our results this is expected, as to find the asymmetric defects one must allow for these small perturbations in far-field spatial wavenumber and wave speed. Building a theoretical framework that accounts for and deals with this additional freedom is the key, novel result of this paper. 

\section{Spatial Dynamics of Wave Trains and Contact Defects} \label{section:spatial_dynamics}

In the previous section, we noted that the essential information concerning snaking bifurcations comes from looking at a reduced 6-dimensional model equation. Of particular importance is the identification of the spatial spectra of the wave trains and Turing patterns from which we derive our geometrical understanding of this problem. In this section, we summarize the rigorous mathematical results underpinning this discussion. We also highlight the relevant properties of the underlying wave trains that will be important in formalising the full asymmetric defect problem.

\subsection{Wave Trains}\label{s:wt}

We start by considering wave train solutions in reaction diffusion systems which represent the asymptotic background states of our defects. Recall our reaction diffusion equation is
\begin{equation}\label{e:rds} 
    u_t = D u_{xx} + f(u).
\end{equation}
\changed{A wave train with spatial wavenumber $k$ and temporal frequency $\omega$ is a solution of \eqref{e:rds} of the form $u(x,t)=u_\mathrm{wt}(kx-\omega t)$ where the associated profile $u_\mathrm{wt}(\phi)$ is a $2\pi$-periodic solution of the equation
\begin{align}\label{eqn:wave train}
    k^2 D \partial_{\phi\phi} u + \omega \partial_\phi u + f(u) = 0. 
\end{align}
Throughout, we assume that $u_\mathrm{wt}^\prime(\phi)\not\equiv0$ is not the zero function, so that the profile depends genuinely on $\phi$.}

\changed{Assume now that $u_\mathrm{wt}^*(k_0x-\omega_0t)$ is a wave train with wavenumber $k_0$ and frequency $\omega_0$. Typically (see below for detailed assumptions), there is then a smooth family $u_\mathrm{wt}(kx-\omega t;k)$ of wave trains of (\ref{e:rds}) that are defined for all $k$ near $k_0$, satisfy $u_\mathrm{wt}(\phi;k_0)=u_\mathrm{wt}^*(\phi)$, and have temporal frequency $\omega$ given by a smooth function $\omega_\mathrm{nl}(k)$ so that $\omega=\omega_\mathrm{nl}(k)$ for each $k$ near $k_0$ with $\omega_0=\omega_\mathrm{nl}(k_0)$. We refer to $\omega_\mathrm{nl}(k)$ as the \emph{nonlinear dispersion relation}. These wave trains travel with phase velocity $c_\mathrm{p}=\omega_\mathrm{nl}(k)/k$ (for $k\neq0$) and transport localized perturbations with the group velocity $c_\mathrm{g}=\partial_k\omega_\mathrm{nl}(k)$ (see \cite{doelman2009} for the latter property).}

\paragraph{Turing patterns:} \changed{We say that a wave train is a Turing pattern if its wavenumber $\kappa\neq0$ is nonzero, its temporal frequency $\omega=0$ vanishes, and its profile $u_\mathrm{Tur}(\phi)$ is even in $\phi$: Turing patterns $u(x,t)=u_\mathrm{Tur}(\kappa x)$ are therefore stationary in time and spatially periodic in $x$. Note that we will use $\kappa$ to denote the spatial wavenumber of Turing patterns. A Turing pattern $u_\mathrm{Tur}(\kappa_0x)$ with wavenumber $\kappa_0\neq0$ is typically part of a family $u_\mathrm{wt}(\kappa x;\kappa)$ of Turing patterns parametrized by their wavenumber $\kappa$ near $\kappa_0$ with $u_\mathrm{Tur}(\phi)=u_\mathrm{wt}(\phi;\kappa_0)$ \citep{devaney1977,vanderbauwhede1992}.}

\paragraph{Spatially homogeneous oscillations:} \changed{We say that a wave train is a spatially homogeneous oscillation (or a homogeneous oscillation for short) if $k=0$ and $\omega=\omega_0\neq0$: these wave train solutions do not depend on $x$, are periodic in $t$ with frequency $\omega_0$, and are therefore of the form $u(x,t)=u_\mathrm{Hom}(-\omega_0 t)$ where the profile $u_\mathrm{Hom}(\phi)$ is a $2\pi$-periodic solution of
\begin{equation}\label{e:homode}
-\omega_0 u_\phi = f(u).
\end{equation}
From now on, we assume that the following hypothesis is met.}

\begin{Hypothesis}[Homogeneous oscillation] \label{hyp:ho}
\changed{Assume that $u_\mathrm{Hom}(-\omega_0 t)$ is a spatially homogeneous oscillation with frequency $\omega_0\neq0$ and that the Floquet multiplier $\rho=1$ of the linearization of \eqref{e:homode} about $u_\mathrm{Hom}$ has algebraic multiplicity one.}
\end{Hypothesis}

\changed{Hypothesis~\ref{hyp:ho} implies that the homogeneous oscillation arises as part of a smooth one-parameter family of wave trains $u_\mathrm{wt}(kx-\omega t;k)$ that is parametrized by their wavenumber $k$ near $0$ (see, for instance, \cite[\S3.3]{sandstede2004}). We focus on the linearization of \eqref{e:rds} about the wave trains $u_\mathrm{wt}(kx-\omega t;k)$ for $k$ close to zero. We emphasize that the results below hold more generally, but we shall need these only for $k$ near $0$ and refer to \cite[\S3.1-3.3]{sandstede2004} for the general case.}

\paragraph{Dispersion relations of wave trains with small wavenumber:}
\changed{The linearization of the reaction-diffusion system \eqref{e:rds} about the given wave train $u(x,t)=u_\mathrm{wt}(kx-\omega t,k)$ with $\omega\neq0$ in the transformed time variable $\tau=\omega t$ is given by
\begin{align}\label{e:lrds}
\omega u_\tau = D \partial_{xx} u + f'(u_\mathrm{wt}(kx-\tau,k))u
\end{align}
which has coefficients that are $2\pi$-periodic in $\tau$. The treatment in \cite[\S3.1-3.3]{sandstede2004} shows that there is a discrete number of analytic curves $\lambda_j(\nu)$ in the complex plane that are each parametrized by $\nu\in\rmi\R$ so that the following is true: $\rho\in\C$ is in the spectrum of the linear time-$2\pi$ map of \eqref{e:lrds} posed on $L^2(\R)$ if and only if $\rho=\exp(2\pi\lambda/\omega_0)$ with $\lambda=\lambda_j(\nu)$ for some $j$ and some $\nu\in\rmi\R$. Moreover, since we assumed that the Floquet multiplier $\rho=1$ of the linearization of \eqref{e:homode} about the homogeneous oscillation is algebraically simple and since $k$ is close to $0$, there is a exactly one curve with $\lambda_j(0)=0$. We refer to the resulting curve, which is defined for $\nu\in\rmi\R$ near $0$, as the \emph{linear dispersion relation} $\lambda_\mathrm{lin}(\nu)$ as it relates the wavenumber and frequency of solutions to the linear system \eqref{e:lrds}. Using Taylor expansion arguments \citep{doelman2009} one can show that the linear and nonlinear dispersion relations are related through the group velocity $c_\mathrm{g}(k)=\partial_k\omega_\mathrm{nl}(k)$ via
\begin{align}
    \omega_\mathrm{nl}(k) &= \omega_0 + c_\mathrm{g}(k_0)(k-k_0) + \bigo(|k-k_0|^2), \label{eqn:nonlinear_dispersion} \\
    \lambda_\mathrm{lin}(\nu;k) &= c_\mathrm{g}(k)\nu + d_\|(k) \nu^2 + \bigo(\nu^3), \label{eqn:linear_dispersion}
\end{align}
for $k$ close to $k_0$ and $\nu$ close to 0. Through \eqref{eqn:linear_dispersion} we see that there is a direct link between the temporal stability properties of a wave train, its group velocity, and its spatial Floquet exponents $\nu$. We summarize our discussion in the following lemma and provide an expansion of the nonlinear dispersion relation at $k=0$.}

\begin{Lemma}[{\cite[\S3.3]{sandstede2004}}] \label{lemma:hom_osc}
\changed{Assume that Hypothesis~\ref{hyp:ho} holds, then there exists a smooth family $u_\mathrm{wt}(kx-\omega_\mathrm{nl}(k)t);k)$ of wave train solutions of \eqref{e:rds} near the homogeneous oscillation $u(x,t)=u_\mathrm{Hom}(-\omega_0t)$ that is parametrized by the wavenumber $k$ for $k$ near $0$ and whose nonlinear dispersion relation $\omega_\mathrm{nl}(k)$ and linear dispersion relation $\lambda_\mathrm{lin}(\nu;0)$ at $k=0$ satisfy
\begin{align}
\omega_\mathrm{nl}(k) = \omega_0 + \frac12 \omega_\mathrm{nl}^{\prime\prime}(0) k^2 + \bigo(k^4), \\
\lambda_\mathrm{lin}(\nu;0) = d_\| \nu^2 + \bigo(\nu^3).
\end{align}
for some smooth function $d_\|(k)$.}
\end{Lemma}

\changed{In particular, we see that the group velocity of spatially homogeneous oscillations vanishes. Finally, we note that we are ultimately interested in solutions to the reaction-diffusion equation with a uniform transport term found by changing variables into a comoving frame,
\begin{equation}
    u_t = Du_{yy} + c_\mathrm{d} u_y + f(u),\label{eqn:rdc_equationy}
\end{equation}
as we will need to allow our defect patterns to move at a non-zero wave speed $c_\mathrm{d}$. Note that a wave train solution $u_\mathrm{wt}(kx-\omega_\mathrm{nl}(k)t;k)$ of \eqref{e:rds} will generate a periodic wave train solution $u(y,t)=u_\mathrm{wt}(ky-\omega_\mathrm{d}t;k):=u_\mathrm{wt}(ky-(\omega_\mathrm{nl}(k)-c_\mathrm{d}k)t;k)$ of the reaction-diffusion equation \eqref{eqn:rdc_equationy} in the comoving frame with the same wavenumber $k$ and the adjusted temporal frequency given by
\begin{equation}\label{e:adjfreq}
\omega_\mathrm{d}(k) := \omega_\mathrm{nl}(k) - c_\mathrm{d}k,
\end{equation}
resulting in the new phase and group velocities given by $\tilde{c}_\mathrm{p}=c_\mathrm{p}-c_\mathrm{d}$ and $\tilde{c}_\mathrm{g}=c_\mathrm{g}-c_\mathrm{d}$, respectively. From now on, we will always use $x$ instead of $y$ and refer to the system
\begin{equation}
u_t = Du_{xx} + c_\mathrm{d} u_x + f(u),\label{eqn:rdc_equation}
\end{equation}
regardless of the frame we are working with.}

\subsection{Spatial Dynamics}

\changed{We assume that Hypothesis~\ref{hyp:ho} holds. We then know that $u_\mathrm{Hom}(-\omega_0t)$ is part of a family $u_\mathrm{wt}(kx-\omega_\mathrm{nl}(k)t);k)$ of wave train solutions of \eqref{e:rds} that is parametrized by the wavenumber $k$ for $k$ near $0$ and satisfies $u_\mathrm{wt}(-\omega_\mathrm{nl}(0)t);0)=u_\mathrm{Hom}(-\omega_0t)$. In a  coordinate frame that moves with a given speed $c_\mathrm{d}$, these wave trains then have temporal frequency $\omega_\mathrm{d}(k)=\omega_\mathrm{nl}(k)-c_\mathrm{d}k$ given by \eqref{e:adjfreq}. From now on, we choose a wavenumber $k$ near zero and a speed $c_\mathrm{d}$ and denote by $\omega_\mathrm{d}$ the resulting adjusted frequency of the wave train with wavenumber $k$, where we do not explicitly label its dependence on $k$.}

So far we have described the spectral properties of time-periodic solutions to the reaction-diffusion equation. In order to leverage dynamical systems arguments, we need access to dynamical information of the elliptic formulation
\begin{align}
    u_x &= w, \label{eqn:rdc_elliptic} \\ \nonumber 
    w_x &= D^{-1}(\omega_\mathrm{d} u_\tau - c_\mathrm{d} w - f(u)),
\end{align}
where we recall that the relevant phase space is $X = H^{1/2}_{\mathrm{per}}( S^1, \R^2)\times L^2_{\mathrm{per}}( S^1, \R^2)$, and we have rescaled time to $\tau = \omega_\mathrm{d} t$ to ensure we have $2\pi$-periodic solutions in $\tau$. 
\changed{The wave trains $u_\mathrm{wt}(kx-\tau;k)$ generate solutions $(u_\mathrm{wt},k u_\mathrm{wt}^\prime)(kx-\cdot;k)\in X$ of \eqref{eqn:rdc_elliptic} which are equilibria in $x$ for $k=0$ and $2\pi/k$-periodic in $x$ for $k\neq0$.} The temporal Floquet problem for the linearization of \eqref{eqn:rdc_elliptic} around the periodic solution $(u_\mathrm{wt},k u_\mathrm{wt}^\prime)(kx-\cdot;k)$ is
\begin{align}
    u_x &= w , \label{eqn:temporal_floquet_problem_1} \\ \nonumber
    w_x &= D^{-1}(\omega_\mathrm{d} u_\tau - c_\mathrm{d} w - f'( u_\mathrm{wt}(kx-\tau;k))u + \lambda u).
\end{align}

Again calling on spatial Floquet theory \citep{sandstede2004, mielke1996}, for every $\lambda \in \C$ the dynamical system \eqref{eqn:temporal_floquet_problem_1} posed on the phase space $X$, has an exponential trichotomy. That is, there exist smooth projection maps $\prwt{u}(x_0, \lambda)$, $\prwt{s}(x_0, \lambda)$ and $\prwt{c}(x_0, \lambda)$ mapping $X\rightarrow X$ such that 
\begin{compactitem}
    \item $\prwt{u}(x_0, \lambda) + \prwt{s}(x_0, \lambda) + \prwt{c}(x_0, \lambda) = I$ for all $(x_0, \lambda)$, 
    \item $\Rg(\prwt{u}(x_0, \lambda))$ is the set of initial conditions at $x = x_0$ for which the dynamical system \eqref{eqn:temporal_floquet_problem_1} converges exponentially to 0 as $x \rightarrow -\infty$, 
    \item $\Rg(\prwt{s}(x_0, \lambda))$ is the set of initial conditions at $x = x_0$ for which the dynamical system \eqref{eqn:temporal_floquet_problem_1} converges exponentially to 0 as $x \rightarrow \infty$, and
    \item $\Rg(\prwt{c}(x_0, \lambda))$ is the set of initial data at $x = x_0$ for which the dynamical system \eqref{eqn:temporal_floquet_problem_1} grow at most algebraically for $x\in\R$.
\end{compactitem}

Moreover, the ranges of $\prwt{u}$ and $\prwt{s}$ are infinite-dimensional, while the range of $\prwt{c}$ is finite-dimensional. Note that when $\lambda = 0$ the temporal Floquet problem reduces to the linearization of the wave train itself. So this establishes the existence of exponential trichotomies for the wave train. In particular, by using exponential weighting arguments as needed, it is possible to construct centre-stable and center-unstable manifolds to the wave train $ u_\mathrm{wt}(kx-\tau;k)$, both of which are infinite-dimensional. The infinite dimension of these manifolds results from the fact that the elliptic problem \eqref{eqn:rdc_elliptic} is ill-posed as an initial value problem on the phase space $X$. In order to recover the finite-dimensional model system described in \S\ref{section:intuition_infinitedim}, we need to be careful. 

To get an accurate picture of the system, we need to understand the structure of the spatial Floquet exponents closest to the imaginary axis. \changed{The equilibria or periodic solutions $(u_\mathrm{wt},k u_\mathrm{wt}^\prime)(kx-\alpha-\cdot;k)\in X$ of \eqref{eqn:rdc_elliptic} belonging to the wave trains come in a one-parameter family parametrized by the phase $\alpha\in S^1$ due to the $S^1$-symmetry provided by translation in the time variable $\tau$. In particular, $(u_\mathrm{wt},k u_\mathrm{wt}^\prime)$ is a relative equilibrium with respect to this symmetry and \eqref{eqn:temporal_floquet_problem_1} with $\lambda=0$ will therefore have a spatial Floquet exponent at zero since the  center manifold will include all time translates $(u_\mathrm{wt},k u_\mathrm{wt}^\prime)(kx-\alpha-\cdot;k)$. The question is then what the next closest spatial Floquet exponents would be. The answer depends on the group velocity $c_\mathrm{g}$ of the wave train with wavenumber $k$ and the chosen speed $c_\mathrm{d}$. The following lemma is formulated again in the context of spatially homogeneous oscillations and the wave trains with wavenumber $k$ close to zero that accompany them though the result is true more generally.}

\begin{Lemma}[{\cite[Lemma~3.6]{sandstede2004}}] \label{lemma:spectra_wave trains}
\changed{Assume that $u_\mathrm{Hom}(-\omega_0t)$ is a spatially homogeneous oscillation for $\omega_0\neq0$ with the following properties:
\begin{compactitem}
\item The Floquet multiplier $\rho=1$ of the linearization of \eqref{e:homode} about $u_\mathrm{Hom}(\phi)$ is algebraically sinmple.
\item The linear dispersion relation $\lambda_\mathrm{lin}(\nu)$ of the homogeneous oscillation satisfies $\lambda_\mathrm{lin}^{\prime\prime}(0)>0$.
\item The wave train is temporally stable for \eqref{e:rds}, that is, with the exception of the Floquet exponents near $\lambda=0$ given by the linear dispersion relation, its Floquet-exponent spectrum is contained in the open left half-plane.
\item The nonlinear dispersion relation is genuinely nonlinear so that $\omega_\mathrm{nl}^{\prime\prime}(0)\neq0$.
\end{compactitem}
Then the spatial dynamical system \eqref{eqn:temporal_floquet_problem_1} belonging to the wave train $u_\mathrm{wt}(kx-\omega t;k)$ for $k$ near zero, which are guaranteed to exist by Lemma~\ref{lemma:hom_osc}, has a geometrically simple Floquet exponent $\nu = 0$ for $\lambda = 0$. This Floquet exponent is algebraically simple if $c_\mathrm{d} \neq c_\mathrm{g}(k)$, while it has algebraic multiplicity two if $c_\mathrm{d} = c_\mathrm{g}(k)$.}
\end{Lemma}

\begin{figure}
    \centering
    \includegraphics[width=0.8\linewidth]{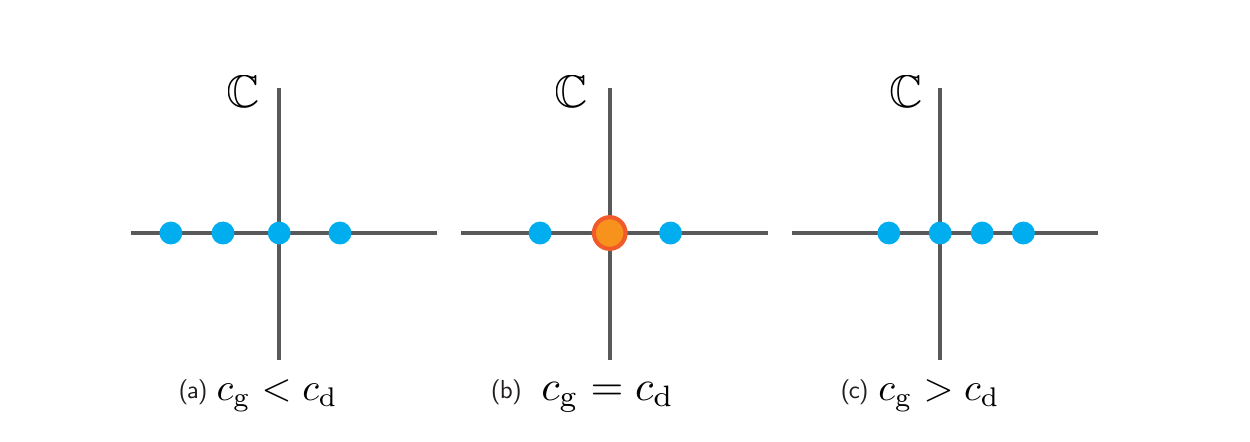}
    \caption[Spatial spectrum on wave trains]{The spatial spectrum of wave trains near the origin for the three cases: (a) $c_\mathrm{g} < c_\mathrm{d}$, (b) $c_\mathrm{g} = c_\mathrm{d}$ and (c) $c_\mathrm{g} > c_\mathrm{d}$.}
    \label{fig:spatial_spectrum_wave trains}
\end{figure}

A direct result of this Lemma is that, near the origin, we have the spatial Floquet spectrum pictures shown in Figure~\ref{fig:spatial_spectrum_wave trains}. With the exception of the first zero Floquet exponent and the second closest Floquet exponent, which crosses through zero as $c_\mathrm{d}-c_\mathrm{g}(k)$ crosses through zero, the remaining Floquet spectra are distributed evenly in the right and left half-planes. They also remain bounded away from the imaginary axis. 

The next lemma states that for each given $c_\mathrm{d}$ near zero there is a unique wavenumber $k$ near zero for which the wave train $u_\mathrm{wt}(kx-\omega_\mathrm{nl}(k)t;k)$ that accompanies a given homogeneous oscillation has group velocity $c_\mathrm{d}$.

\begin{Lemma} \label{lemma:group_velocities}
Assume that the hypotheses of Lemma~\ref{lemma:spectra_wave trains} hold for the homogeneous oscillation $u_\mathrm{Hom}(-\omega_0t)$. Then for every $c$ near zero there exists a unique wavenumber $k$ near zero such that the wave train $u_\mathrm{wt}(kx-\omega_\mathrm{nl}(k)t; k)$ guaranteed by Lemma~\ref{lemma:hom_osc} has group velocity $c$. Moreover, we have the leading-order expansions
\begin{align}
k &= \frac{c}{\omega_\mathrm{nl}^{\prime\prime}(0)} + \bigo(c^3), \\
\omega(c) &:= \omega_\mathrm{nl}(k(c)) = \omega_0 + \frac{c^2}{2\omega_\mathrm{nl}^{\prime\prime}(0)} + \bigo(c^4).
\end{align}
\end{Lemma}

\begin{Proof}
From Lemma~\ref{lemma:hom_osc}, we know that the nonlinear dispersion relation of the family of wave trains generated near the homogeneous oscillation is given by
\begin{align}
\omega_\mathrm{nl}(k) = \omega_0 + \frac12 \omega_\mathrm{nl}^{\prime\prime}(0) k^2 + \bigo(k^4),
\end{align}
and differentiating with respect to $k$ we obtain the group velocity
\begin{equation}
c_\mathrm{g}(k) = \frac{\mathrm{d}\omega_\mathrm{nl}}{\mathrm{d}k}(k) = \omega_\mathrm{nl}^{\prime\prime}(0) k + \bigo(k^3).
\end{equation}
We solve the equation $c_\mathrm{g}(k) = c$ for $k$ as a function of $c$ using the Implicit Function Theorem which also yields the expansions stated above.
\end{Proof}

\changed{In particular, for each $c_\mathrm{d}$ near zero there is a unique wavenumber $k$ near zero so that the center manifold of the periodic solution $(u_\mathrm{wt},ku_\mathrm{wt}^\prime)(kx-\cdot;k)$ of \eqref{eqn:rdc_elliptic} has dimension two. The latter is the crucial property that will allow us to construct symmetric and asymmetric contact defects, as the additional center-direction produced by the group-velocity condition is necessary for their formation.}

\subsection{Geometry and Robustness of Contact Defects}\label{s:geom}

\begin{figure}
    \centering
    \includegraphics[width=0.7\textwidth]{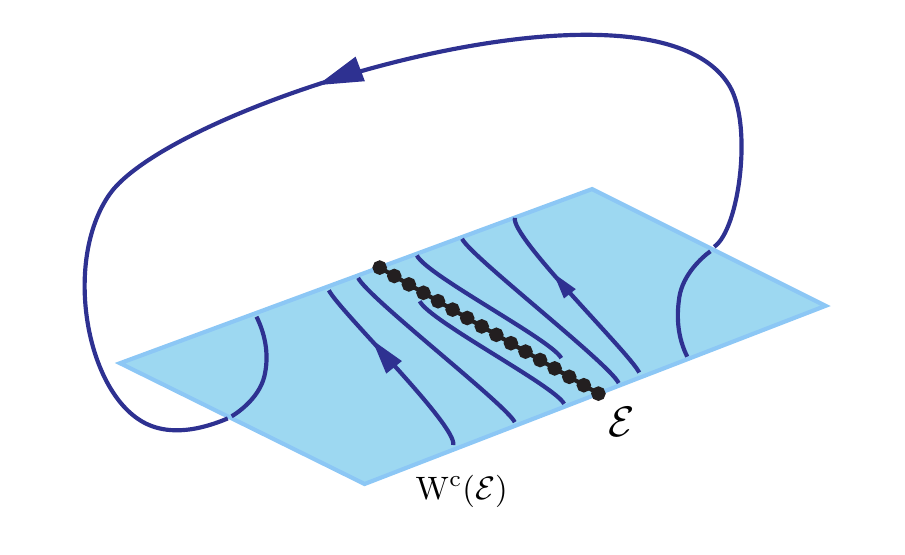}
    \caption[Geometry of contact defects]{The geometry of contact defects in the phase space $X$. Contact defects are homoclinic connections to the invariant circle of relative equilibria made through the two-dimensional center manifold.}
    \label{fig:geometry_contact_defects}
\end{figure}

\changed{We first discuss the definition and geometric properties of standing contact defects. Assume that the hypotheses of Lemma~\ref{lemma:spectra_wave trains} hold for the homogeneous oscillation $u_\mathrm{Hom}(-\omega_0t)$ of \eqref{eqn:rdc_equation} with $c_\mathrm{d}=0$. Note that Lemma~\ref{lemma:hom_osc} shows that these oscillations have group velocity zero, and Lemma~\ref{lemma:spectra_wave trains} therefore implies that the circle $\mathcal{E}:=(u_\mathrm{Hom}(\cdot-\alpha),0)_{\alpha\in S^1}\in X$ of equilibria of the spatial dynamical system \eqref{eqn:rdc_elliptic} with $c_\mathrm{d}=0$ and $\omega_\mathrm{d}:=\omega_0$ has a two-dimensional center manifold. As shown in \cite{doelman2009, sandstede2004a}, the flow on the center manifold is given by
\begin{equation}\label{e:vfcm}
\theta_x = k, \qquad
k_x = -\frac{k^2}{a} + g(k) k^4, \qquad
a := -\frac{\lambda_\mathrm{lin}^{\prime\prime}(0)}{\omega_\mathrm{nl}^{\prime\prime}(0)} \neq 0
\end{equation}
for some smooth function $g(k)$, where $\theta\in S^1$ corresponds to the phase of the homogeneous oscillation, and $k$ spans the remaining direction in the center manifold and corresponds to the spatial wavenumber. Thus, for $a>0$, solutions in the center manifold with $k>0$ converge to the circle $\mathcal{E}=\{k=0\}$ of equilibria as $x\to\infty$, while solutions with $k<0$ converge to this circle as $k\to\infty$; see Figure~\ref{fig:geometry_contact_defects} for an illustration. For $a<0$, the same statement is true in opposite directions of $x$. Standing contact defects are solutions $(u_\mathrm{d},\partial_x u_\mathrm{d})(x,\cdot)\in X$ of \eqref{eqn:rdc_elliptic} that lie in the intersection of the center-unstable manifold $W^\mathrm{cu}(\mathcal{E})$ and the center-stable manifold $W^\mathrm{cs}(\mathcal{E})$ of the circle $\mathcal{E}$ of equilibria and that converge to this circle as $x\to\pm\infty$. They can therefore be thought of as homoclinic orbits to the circle $\mathcal{E}$. Since the two center directions are counted towards both $W^\mathrm{cu}(\mathcal{E})$ and $W^\mathrm{cs}(\mathcal{E})$, we expect that they generically intersect transversely along a two-dimensional manifold, which is then given by the space and time translates of the contact defect $(u_\mathrm{d},\partial_x u_\mathrm{d})(x,\cdot)$. This statement was made rigorous in \cite[Lemma~5.6]{sandstede2004}. As a consequence, this intersection is robust and persists under changes of a system parameter. It is also robust when we replace the circle of equilibria consisting of spatially homogeneous oscillations by the periodic orbit $\mathcal{P}(k):=\{ (u_\mathrm{wt},k u_\mathrm{wt}^\prime)(kx-\cdot;k)\colon x\in\R\}$ for a wavenumber $k$ close to zero provided we maintain the saddle-node structure of the flow on the center manifold. Given $k$ near zero, we therefore need to select the defect speed $c_\mathrm{d}$ to coincide with the group velocity of the selected wave train and adjust the temporal frequency $\omega_\mathrm{d}$ via \eqref{e:adjfreq} to account for the moving frame by setting
\[
c_\mathrm{d} := \omega_\mathrm{nl}^\prime(k) = \omega_\mathrm{nl}^{\prime\prime}(0) k + \bigo(k^3), \qquad
\omega_\mathrm{d} := \omega_\mathrm{nl}(k) - c_\mathrm{d}k = \omega_\mathrm{nl}(k) - \omega_\mathrm{nl}^\prime(k) k.
\]
The resulting contact defects travel with speed $c_\mathrm{d}$ and converge to the periodic orbit $\mathcal{P}(k)$ as $x\to\pm\infty$.}

\changed{We now discuss the convergence properties of the standing contact defects towards the homogeneous oscillations in the far field as $|x|\to\infty$. In the spatial-dynamics formulation, the contact defects $(u_\mathrm{d},\partial_x u_\mathrm{d})(x,\cdot)$ approach the circle $\mathcal{E}$ of equilibria, but they will not converge towards a specific equilibrium, that is, they will not converge with asymptotic phase. To see this, we solve the vector field \eqref{e:vfcm} on the center manifold to get
\begin{equation}\label{etheta}
\theta_\pm(x) = \theta_\pm^0 + a \log |x| + \bigo(1/|x|)
\end{equation}
as $x\to\pm\infty$, where the logarithmic term is universal regardless of whether $x\to\infty$ or $x\to-\infty$, while the phases $\theta_\pm^0$ may depend on the sign of $x$. As already shown in \cite{sandstede2004a}, this implies the asymptotics
\[
|u_\mathrm{d}(x,t)-u_\mathrm{Hom}(-\omega_0 t+\theta_\pm(x))| \rightarrow 0 \quad\mbox{ as }\quad x\to\pm\infty
\]
for the contact defect $u_\mathrm{d}(x,t)$ of \eqref{e:rds}. The logarithmic behavior of the phase function is an important consideration when constructing numerical schemes to approximate contact defects.}

\changed{We mention two consequences of the discussion above. Firstly, \eqref{etheta} shows that the phase difference $\theta_+^0-\theta_-^0$ of the homogeneous oscillations in the far field across the core region is well defined. In the parlance of spiral and target waves, $\theta_+^0-\theta_-^0=0$ modulus $2\pi$ gives a one-dimensional target pattern, while $\theta_+^0-\theta_-^0=\pi$ modulus $2\pi$ gives a one-dimensional spiral wave. Secondly, we see that the slopes of the curves $-\omega_0t+\theta_\pm(x)=\mathrm{const.}$ of constant phase of the contact-defect profile in the far field of the space-time plane are given by
\begin{equation}\label{e:velocity}
\frac{\mathrm{d}t}{\mathrm{d}x} = \frac{a}{\omega_0 x} + \bigo(1/x^2), \qquad |x|\gg1
\end{equation}
instead of decaying exponentially as for source defects.}

\subsection{Turing Patterns}

\changed{Finally we consider the spatial-dynamics properties of Turing patterns. Recall from \S\ref{s:wt} that Turing patterns are stationary, spatially periodic solutions $u(x,t)=u_\mathrm{Tur}(\kappa x;\kappa)$ of \eqref{e:rds} that arise generically in one-parameter families parametrized by the nonzero wavenumber $\kappa\neq0$, where $u_\mathrm{Tur}(\phi;\kappa)$ is $2\pi$-periodic and even in $\phi$ for each wavenumber $\kappa$ \citep{devaney1977,vanderbauwhede1992}. Hence, Turing patterns correspond to $2\pi/\kappa$-periodic solutions $(u_\mathrm{Tur},\kappa u_\mathrm{Tur}^\prime)(x;\kappa)$ of \eqref{eqn:rdc_elliptic}. This family forms a two-dimensional center manifold of \eqref{eqn:rdc_elliptic} that is parametrized by the spatial evolution variable $x$ and the wavenumber $\kappa$. The linearization of \eqref{eqn:rdc_equation} with wave speed $c_\mathrm{d}=0$ and any chosen temporal frequency $\omega_\mathrm{d}\neq0$ is given by}
\begin{align}
u_x &= w ,\\
w_x &= D^{-1}(\omega_\mathrm{d} u_\tau - f'(u_\mathrm{Tur}(x;\kappa))u).
\end{align}
This system has a spatial Floquet exponent zero with geometric multiplicity one (with eigenfunction $(u_\mathrm{Tur}', \kappa u_\mathrm{Tur}^{\prime\prime})$) and algebraic multiplicity two, reflecting the parametrization of the two-dimensional center manifold by $\kappa$. We can determine the remaining Floquet exponents by using a Floquet ansatz, $(u,w) = \mathrm{e}^{\nu x}(\bar{u},\bar{w})$ for $2\pi/\kappa$-periodic functions $\bar{u}$ and $\bar{w}$. That is, we want to solve
\begin{align}
    \bar{u}_x &= \bar{w} - \nu \bar{u}, \\
    \bar{w}_x &= D^{-1}(\omega_\mathrm{d}\partial_\tau - f'(u_\mathrm{Tur}(x;\kappa)))\bar{u} - \nu \bar{w}
\end{align}
for $(\bar{u}, \bar{w})$ as $2\pi/\kappa$-periodic functions. We can simplify these equations by taking a temporal Fourier transform which decouples the system into countably many, four-dimensional subproblems for the Fourier coefficients of each order. The identified center Floquet multipliers lie in the temporally homogeneous subspace.

When we move into the comoving frame, things look different. While it is true that we can find solutions of the form $u(x,\tau) = u_\mathrm{Tur}(x+c\tau/\omega_\mathrm{d}; \kappa)$ of the reaction diffusion equation \eqref{eqn:rdc_equation}, they are no longer $2\pi$-periodic in $\tau$ and hence do not lie in the phase space $X$. However we note that the two-dimensional center manifold is normally hyperbolic when $c_\mathrm{d} = 0$, and so will perturb to a new invariant two-dimensional manifold for small $c_\mathrm{d}$, whose dynamics we can express using a Fenichel normal form. The resulting manifold will no longer be made up of periodic orbits, but will generate a modulational drift in the Turing wavenumber that will allow us to produce asymmetric defects.

\section{Set up and Hypotheses} \label{section:set_up_and_hyp}

We begin by setting the ground for the proofs of Theorems~\ref{thm:symmetric_contact} and~\ref{thm:asymmetric_contact}. These can broadly be divided into two groups: local assumptions on the Turing family and background states; and global assumptions on how the background states interact with a neighborhood of the Turing family via their invariant manifolds. For the remainder of this section we will simplify our notation so that \eqref{eqn:spatial_dynamics} is written as
\begin{align}
    \mathbf{u}_x = F(\mathbf{u};\mu,c,\omega) =: A(\mu,c,\omega) \mathbf{u} + N(\mathbf{u}; \mu) \label{eqn:ode}
\end{align}
where $\mathbf{u} = (u,w)^T$, $A$ is the linear part, $N$ the remaining nonlinear terms and $\mu$ is our primary bifurcation parameter. The system parameters $c$ and $\omega$ will be chosen in neighborhoods of $0$ and $\omega_0$ respectively. We will consider $\mu$ to take values in a interval of $\R$ with non-empty interior $J$.

\begin{Hypothesis}[Exponential Dichotomies] \label{hyp:resolvent}
    We assume there exists a constant $C$ such that
    \begin{align}
        \|(\rmi k - A(\mu, c, \omega))^{-1}\|_{L(X)} \leq \frac{C}{|k| + 1 }
    \end{align}
    for all $k \in \R$. We also assume that there exists a $\delta > 0$ such that $|\Re \lambda | > \delta$ for all $\lambda \in \sigma(A)$.
\end{Hypothesis}

Hypothesis~\ref{hyp:resolvent} is a technical assumption required to guarantee the existence, continuation and perturbation of exponential dichotomies \citep{peterhof1997} as described in Section~\ref{section:spatial_dynamics}. All our following results, including existence of invariant manifolds, Fenichel normal form and exchange lemma, require the presence of exponential dichotomies and hence this assumption. We will not directly use this in our analysis.

To begin in earnest we note that our system \eqref{eqn:spatial_dynamics} is spatially-reversible and temporally-equivariant under two distinct symmetry operations which will play key roles for our framework, so we give precise definitions below.

\begin{Definition}[Reversibility and Equivariance] \label{hyp:reversible}
    The dynamical system \eqref{eqn:ode} is {\bf reversible} if there exists an involution $\mathcal{R} : X \rightarrow X$ such that $F(\mathcal{R}\mathbf{u};\mu) = - \mathcal{R} F(\mathbf{u}; \mu)$ for all $\mathbf{u}\in X$. The dynamical system $\eqref{eqn:ode}$ is {\bf equivariant} with respect to a family of transformations $\Gamma_\alpha:X\rightarrow X$ for $\alpha\in S^1$ if $F(\Gamma_\alpha \mathbf{u};\mu) = \Gamma_\alpha F(\mathbf{u};\mu)$ for all $\alpha$ in $S^1$ and $\Gamma_\alpha \circ \Gamma_{\tilde{\alpha}} = \Gamma_{\alpha+\tilde{\alpha}}$ for all $\alpha,\tilde{\alpha}$ in $S^1$.
\end{Definition}

As Equation \eqref{eqn:spatial_dynamics} is a reaction diffusion equation written as a first order spatial dynamical system, it is reversible with $\mathcal{R}$ defined by $(u,w)\mapsto(u,-w)$ and equivariant with respect to the temporal shift operator $[\Gamma_\alpha \mathbf{u}](t) := \mathbf{u}(t+\alpha)$. The combination of these symmetries also generates a second reverser for equation \eqref{eqn:ode}; namely the map $\mathcal{R}\Gamma_\pi = \Gamma_\pi \mathcal{R}$ satisfies the above condition for reversibility. This second reverser will play an important role in constructing multiple symmetric defect branches.

\begin{Hypothesis}[Background Oscillations] \label{hyp:background}
\changed{For each $\mu \in J$ there exists a circle $\mathcal{E}(\mu) = (\Gamma_\alpha\mathbf{u}_\mathrm{Hom})_{\alpha\in S^1}=(u_\mathrm{Hom}(-\omega_0(\mu)t+\alpha;\mu),0)_{\alpha\in S^1}$ of homogeneous oscillations with temporal period $\omega_0(\mu)\neq0$. These homogeneous oscillations satisfy the assumptions of Lemma~\ref{lemma:spectra_wave trains} so that they have a spatial Floquet spectrum with two eigenvalues at zero followed by $\pm\lambda_1(\mu)$ and $\pm \lambda_2(\mu)$ with $0<\lambda_1(\mu),\lambda_2(\mu)$ uniformly in $\mu$ and all other spectrum bounded uniformly outside of this spectrum. The homogeneous oscillations are also accompanied by families of wave trains $u_\mathrm{wt}(kx-\omega(\mu)t;k,\mu)$ with frequency $\omega(\mu)=\omega_\mathrm{nl}(k,\mu)=\omega_0(\mu)+\frac12\omega_\mathrm{nl}^{\prime\prime}(0,\mu)k^2+\bigo(k^4)$.}
\end{Hypothesis}

By equivariance, translation in time $\Gamma_\alpha$ parameterizes the invariant circle \changed{$\mathcal{E}(\mu)$} in phase space:
\begin{align}
    (u_\mathrm{Hom},0)(\cdot+\alpha;\mu) = [\Gamma_\alpha (u_\mathrm{Hom},0)](\cdot;\mu).
\end{align}
Throughout our discussion, we must be careful to keep track of the symmetry implications of this equivariance.

\begin{Hypothesis}[Turing Family] \label{hyp:Turing}
    For every $\mu \in J$ and $\kappa \in K \subset \R^+$ there exists a family of periodic orbits $\mathbf{u}_\mathrm{Tur}(x;\kappa,\mu)$ smooth in $(\kappa, \mu)$ and $\mathcal{R}$-reversible. The weakest spatial Floquet spectrum of these states consist of two zero Floquet exponents, and two pairs of complex conjugate exponents $\pm \rho(\kappa,\mu) \pm \rmi\beta(\kappa,\mu)$ with real parts bounded uniformly away from zero. All other spectrum is bounded uniformly away from this spectrum. Moreover, we assume the eigenspaces corresponding to the weak stable and unstable exponents only intersect the space of time homogeneous functions trivially at zero.
\end{Hypothesis}

We get a 1-parameter family due to reversibility as noted in Section~\ref{section:spatial_dynamics}. The final assumption about the eigenspaces is also essential to the following analysis as it allows us to prescribe how the equivariance $\Gamma_\alpha$ operates on the stable and unstable fibers of the Turing family. Recall from Section~\ref{section:spatial_dynamics}, the Floquet eigenspaces can be identified with the temporal Fourier modes of solutions to the linearized problem. The system is decoupled so our assumption on eigenspaces is that the weakest directions do not correspond to time independent functions as is the case in the Brusselator. 

Importantly for us this has implications on how the equivariance behaves in a neighborhood of the Turing patterns, even though these patterns themselves are temporally homogeneous. Let the eigenfunction corresponding to the weakest stable Floquet exponent $\nu = -\rho + \rmi\beta$ be given by $(\bar{u},\bar{w})$ (and its complex conjugate). By our assumption, these functions are not constant and by the set up described above we can express this function as $\zeta_0(x) \mathrm{e}^{(-\rho+\rmi\beta)x} \mathrm{e}^{2\pi n \rmi t}$ for some $n$ where $\zeta_0$ is a spatially periodic function. We can act on this eigenfunction by our equivariance operator $\Gamma_\alpha$. Time dependence is captured entirely by the exponential prefactor $\mathrm{e}^{2\pi \rmi n t}$ meaning $\Gamma_\alpha$ is a rigid rotation of period $2\pi/n$ in the stable subspace. In particular, this means that we can exactly prescribe the action of $\Gamma_\alpha$ as a $2\pi$-periodic rigid rotation in the weakest stable subspace (and similarly in the weakest unstable subspace). Note that in principle we do not make any assumption on $2\pi$ being the minimal period of the rotation, only that it is a period. Without loss of generality, we will assume that the weakest direction is in the $n=1$ subspace and consequently that the period of $\Gamma_\alpha$ is $2\pi$.

With these local assumptions on the Turing family we can construct a Fenichel normal form and exchange lemma which will provide the framework for our proof.

\begin{Lemma}[Fenichel Normal Form] \label{lemma:fenichel}
    Assume Hypotheses~\ref{hyp:resolvent},~\ref{hyp:reversible} and~\ref{hyp:Turing}. In a neighborhood $\mathcal{V}$ of the Turing family $u_\mathrm{Tur}(\cdot; \cdot, \mu)$, there exists a smooth, $\mathcal{R}$-reversible and equivariant change of coordinates, uniform in $\omega$ near $\omega_0(\mu)$, such that the leading order dynamics are given by
    \begin{align}
    \begin{split}
        \dot{\varphi} &= \kappa, \\
        \dot{\kappa} &= c, \\
        \dot{v}^\mathrm{s} &= \begin{pmatrix} -\rho & -\beta \\ \beta & -\rho \end{pmatrix} v^\mathrm{s}, \\
        \dot{v}^\mathrm{ss} &= A^\mathrm{ss} v^\mathrm{ss}, \\
        \dot{v}^\mathrm{u} &= \begin{pmatrix} \rho & \beta \\ -\beta & \rho \end{pmatrix} v^\mathrm{u}, \\
        \dot{v}^\mathrm{uu} &= A^\mathrm{uu} v^\mathrm{uu}, \\
    \end{split}
    \end{align}
    where $\varphi, \kappa \in \R$ are the center coordinates (periodic phase and spatial wavenumber), $v^\mathrm{s}, v^\mathrm{u} \in \{ z\in\R^2 : |z|< \delta \}$ are the weakest stable and unstable directions and $v^\mathrm{ss}, v^\mathrm{uu} \in Y$ are the remaining strong stable and unstable directions. 
    Moreover, the reverser is given by
    \begin{align}
        \mathcal{R}(\varphi, \kappa, v^\mathrm{s}, v^\mathrm{ss}, v^\mathrm{u}, v^\mathrm{uu}) = (-\varphi, \kappa, v^\mathrm{u}, v^\mathrm{uu}, v^\mathrm{s}, v^\mathrm{ss}) \label{eqn:ch3_fenichel_rev}
    \end{align}
    and the time-shift is given by
    \begin{align}
        \Gamma_\alpha (\varphi, \kappa, v^\mathrm{s}, v^\mathrm{ss}, v^\mathrm{u}, v^\mathrm{uu}) = (\varphi, \kappa, R_\alpha v^\mathrm{s}, T_\alpha v^\mathrm{ss}, R_\alpha v^\mathrm{u}, T_\alpha v^\mathrm{uu}). \label{eqn:ch3_fenichel_equi}
    \end{align}
    for all $\alpha \in S^1$ where $R_\alpha$ is the rotation matrix through an angle $\alpha$ in $\R^2$ and $T_\alpha$ is a $2\pi$-periodic operator on the strong stable and unstable directions.
\end{Lemma}

The proof of this statement can be found in \cite{beck2009}, where we note the form of the spatial Floquet exponents described above mean we can enforce the actions of $\mathcal{R}$ and $\Gamma_\alpha$ described in \eqref{eqn:ch3_fenichel_rev} and \eqref{eqn:ch3_fenichel_equi}. For $c=0$ the first two equations are the center-manifold reduction of the Turing family where we see that the Jordan block is expressed as a change in the spatial wavenumber of the oscillations in the family. When $c \neq 0$ we use that the center-manifold is normally hyperbolic to find a new 2-dimensional invariant manifold whose dynamics are a regular perturbation $\bigo(c)$-close to the original center dynamics. A Melnikov integral computation shows that for $c$ close to zero, the wavenumber $\kappa$ is adjusted directly by $c$. \changed{ As discussed in Section~\ref{s:geom}, when changing into a comoving frame, we need to suitably adjust the temporal wavenumber to maintain contact defects. To do this, we note that the Fenichel normal form of Lemma~\ref{lemma:fenichel} can be found uniformly in $\omega$ near $\omega_0$. By Lemma~\ref{lemma:group_velocities}, we can then adjust directly to the correct temporal frequency which is $\omega_d(\mu) = \omega_0(\mu) + \bigo(c^2)$ removing direct dependence on $\omega$ in exchange for higher order dependence on $c$.}

It is important to note that this system is legitimately infinite-dimensional. While we have written these coordinates in a form that exposes the six weakest directions, we cannot ignore the higher order terms. What we will show in the following analysis is that, under subsequent assumptions on the transversality and general position of the incoming manifolds from the homoegeneous oscillations, the infinite-dimensional higher order modes $v^\mathrm{ss}$ and $v^\mathrm{uu}$ can be solved for by the Implicit Function Theorem. 

In this coordinate system the Turing family itself is the set $\{ v^\mathrm{s} = v^\mathrm{u} = v^\mathrm{ss} = v^\mathrm{uu} = 0\}$ with $\varphi$ the direction of motion around the Turing family and $\kappa$ the spatial wavenumber. When considering defect solutions that do not transport in the laboratory frame, $\dot{\kappa} = 0$ indicating that there are no dynamics joining Turing patterns of different wavenumber. However, if we move into a comoving frame with some non-zero constant speed $c$ then the invariant manifold becomes dynamic with a monotone increase or decrease in wavenumber over time. In order to produce symmetric defects it suffices to keep $c = 0$. However, asymmetric defects join fronts and backs associated with Turing patterns of different wavenumber, necessitating a nonzero wave speed $c$.

The center-stable and center-unstable manifolds of the Turing family are given by the sets $\{ v^\mathrm{u} = v^\mathrm{uu} = 0\}$ and $\{ v^\mathrm{s} = v^\mathrm{ss} = 0 \}$ respectively. Also the two subspaces
\begin{align}
\begin{split}
    \Sigma^\mathrm{s} &= \{ |v^\mathrm{s}| = \delta \}, \\
    \Sigma^\mathrm{u} &= \{ |v^\mathrm{u}| = \delta \},
\end{split}
\end{align}
form transverse sections to the flow of system \eqref{eqn:ode} and will be critical in our construction of defect solutions. As a result of Lemma~\ref{lemma:fenichel} we can construct a family of solutions joining these two sections which stay close to these center-stable and center-unstable manifolds for large times $L$.

\begin{Lemma}[Exchange Lemma] \label{lemma:exchange}
    Assume Hypotheses~\ref{hyp:resolvent},~\ref{hyp:reversible} and~\ref{hyp:Turing}. Then for every $\mu\in J$, $\kappa_0\in K$, $c$ small enough, angles $\varphi_0, \theta^\mathrm{s}_0, \theta^\mathrm{u}_0 \in S^1$, $v^\mathrm{ss}_0,v^\mathrm{uu}_0 \in Y$ and $L>0$ sufficiently large there is a unique solution $v(x) = v(x; \mu, c, \kappa_0, \varphi_0, \theta^\mathrm{s}_0, \theta^\mathrm{u}_0, v^\mathrm{ss}_0, v^\mathrm{uu}_0)$ such that
    \begin{align}
    \begin{split}
        &v(x) \in \mathcal{V} \quad \text{for all} \quad x\in[-L, L], \\
        &v(-L) \in \Sigma^\mathrm{s}, \quad \text{and} \\
        &v(L)\in \Sigma^\mathrm{u}.
    \end{split}
	\end{align}
    and with the asymptotics
    \begin{align}
        v(-L) &= (\varphi_0-\kappa_0 L+\frac{1}{2}cL^2 + \mathcal{O}, \kappa_0-cL+\mathcal{O}, \delta \hat{e}(\theta^\mathrm{s}_0), v^\mathrm{ss}_0 + \mathcal{O}, \mathcal{O}, \mathcal{O}), \label{eqn:exchange1} \\
        v(L) &= (\varphi_0 + \kappa_0 L + \frac{1}{2}cL^2 + \mathcal{O}, \kappa_0 + cL + \mathcal{O}, \mathcal{O}, \mathcal{O}, \delta \hat{e}(\theta^\mathrm{u}_0), v^\mathrm{uu}_0 + \mathcal{O}), \label{eqn:exchange2}
    \end{align}
    where $\hat{e}(\Theta)$ is the unit vector in $\R^2$ making angle $\Theta$ with the first component axis and $\mathcal{O} = \mathcal{O}(\mathrm{e}^{-\eta L})$ and all such error estimates are smooth and differentiable.
\end{Lemma}

This is a central result stemming from geometric singular perturbation theory \citep{jones1996}, and extended more generally to any normally hyperbolic sets \citep{schecter2008,schecter2008a}. For a proof of this result in the context of spatial dynamics on infinite-cylinders see \cite{beck2009}.

Now we need to make global assumptions on the intersection of the stable and unstable manifolds of the background states with the sections $\Sigma^\mathrm{s}$ and $\Sigma^\mathrm{u}$. From this we can shoot from the unstable manifold to the stable manifold by way of our exchange lemma and enact a Lin's method matching argument to prove the existence of our defect solutions.

\begin{Hypothesis}[Fronts and Invariant Manifolds] \label{hyp:global}
    For $c=0$ and $\omega = \omega_0(\mu)$, there exists a bifurcation family of heteroclinic connections between the circle $\mathcal{E}(\mu)$ of homogeneous oscillations $\mathbf{u}_\mathrm{Hom}(\cdot; \mu)$ and the Turing patterns $\mathbf{u}_\mathrm{Tur}(\cdot; \cdot, \mu)$ as $\mu$ is varied through $J$. By equivariance, for each fixed $\mu$ there is an invariant circle of such connections obtainable by the action of $\Gamma_\alpha$ on any single member of the family.

    Moreover, we assume that there is a smooth periodic function $(z,\kappa^*,v^\mathrm{ss}_*): S^1\to J\times K\times Y$,  $\varphi\mapsto(z,\kappa^*,v^\mathrm{ss}_*)(\varphi)$ so that $(\varphi, \kappa, v^\mathrm{ss},\mu)\in \mathcal{B}$ with
    \changed{\begin{align}
        \mathcal{B} := \{ (\varphi,\kappa, v^\mathrm{ss}, \mu) \ : \ \rmW^\mathrm{cu}(\mathcal{E}(\mu)) \cap \rmW^\mathrm{cs}(\{\mathbf{u}_\mathrm{Tur}(x;\kappa,\mu)\colon x\in\R\})\cap \Sigma^\mathrm{s} \neq \emptyset \}
    \end{align}}
    if, and only if, $\mu = z(\varphi)$, $\kappa = \kappa^*(\varphi)$ and $v^\mathrm{ss} = T_\alpha v^\mathrm{ss}_*(\varphi)$ for some $\alpha \in S^1$.

    Finally we also assume that the intersection of the incoming manifold \changed{$\rmW^\mathrm{cu}(\mathcal{E}(\mu))$} with the section $\Sigma^\mathrm{s}$ can be parameterized as a graph $g$ over $(\varphi, \kappa, \theta^\mathrm{s}, v^\mathrm{uu}, \mu, c)$. That is, \changed{$p\in \rmW^\mathrm{cu}(\mathcal{E}(\mu))\cap \Sigma^\mathrm{s}$} nearby the fronts if, and only if, 
    \begin{align} 
        p = (\varphi, \kappa, \delta \hat{e}(\theta^\mathrm{s}), g^\mathrm{ss}(q), g^\mathrm{u}(q), v^\mathrm{uu}) \label{eqn:pW-}
    \end{align}
    where $q = (\varphi, \kappa,\theta^\mathrm{s}, v^\mathrm{uu}, \mu, c)$. We assume that the function $g$ and its first derivatives are uniformly bounded in $\mu$.
\end{Hypothesis}

The function $z$ encodes the bifurcation diagram for the fronts. The fact that it can be written as a periodic function is the feature that determines the snaking of defect solutions. The existence proofs that follow can be weakened to show that the presence of a front necessitates nearby defect solutions. However, the continuation of these solutions and the qualitative features of the snaking diagram require the stronger periodic assumption.

The function $g = (g^\mathrm{u}, g^\mathrm{ss})$ describes the expected behavior of the incoming and outgoing stable and unstable manifolds. Observe that $g^\mathrm{u}$ takes values in $\R^2$ (isomorphic to the unstable manifold of the Turing family) and $g^\mathrm{ss}$ values in $Y$ (isomorphic to the strong stable fibers of the Turing family). There are many various assumptions one can make here in which the function $g$, and its domain and range, have a different graph structure. This particular choice of chart for the manifolds should remain valid for all $\mu$ in the bifurcation interval. However, when solving the matching equations it requires a slightly stronger condition on the derivatives of $g$ (see Hypothesis~\ref{hyp:degeneracy} below) than is necessary. This allows us to use a single setup to globally construct the snaking diagram. It is possible that such a parameterization would fail at some values of $\mu$, for example, the fast stable and unstable directions may become involved in the bifurcation. We should think of this as non-generic behavior as it would signal the presence of secondary, higher-codimension bifurcations such as inclination flips or resonances in the defect solutions (viewed as homoclinic orbits). We do not consider such cases here. Our assumed parameterization would also fail in the case of a conservative equation as the wavenumber $\kappa$ would be constant (similarly to the Swift--Hohenberg equation). In such cases, it is possible to amend our current set up using other parameterizations which can be employed separately, or patched together with the one described here, with little change to the overall argument. So we assume that we are in a region where the global parameterization holds.

Hypothesis~\ref{hyp:global} does not directly assume any form on the incoming center-stable manifold as we can use reversibility of the Fenichel normal form to give an explicit representation of this manifold in terms of the incoming center-unstable manifold. From Hypothesis~\ref{hyp:global} we know that a point $p$ lies in \changed{$\Sigma^\mathrm{s} \cap \rmW^\mathrm{cu}(\mathcal{E}(\mu))$} nearby a front if, and only if, there exist constants $q_- = (\varphi_-, \kappa_-, \theta^\mathrm{s}, v_-^\mathrm{uu})$ such that
\begin{align}
    p = (\varphi_-, \kappa_-, \delta \hat{e}(\theta^\mathrm{s}), g^\mathrm{ss}(q_-, \mu, c), g^\mathrm{u}(q_-, \mu, c), v_-^\mathrm{uu}).
\end{align}
When we apply the reverser to this point, we must then have that $\mathcal{R}p$ lies in \changed{$\Sigma^{\mathrm{u}}\cap \rmW^\mathrm{cs}(\mathcal{E}(\mu))$} as $\mathcal{R}$ maps \changed{$\Sigma^\mathrm{s}$ to $\Sigma^\mathrm{u}$} and the center-unstable manifold to the center-stable (and vice versa). Computing this action we see that
\begin{align}
    \mathcal{R}p = (-\varphi_-, \kappa_-, g^\mathrm{u}(q_-, \mu, c), v^\mathrm{uu}_-, \delta \hat{e}(\theta^\mathrm{u}), g^\mathrm{ss}(q_-, \mu, c)).
\end{align}

We can then exactly characterize the opposite intersection as \changed{$\bar{p} \in \Sigma^\mathrm{u}\cap W^\mathrm{cs}(\mathcal{E}(\mu))$} nearby a back if, and only if, there exist constants $q_+ = (\varphi_+, \kappa_+, \theta^\mathrm{u}, v_-^\mathrm{ss})$ such that
\begin{align}
    \bar{p} = (-\varphi_+, \kappa_+, g^\mathrm{u}(q_+, \mu, c), v_+^\mathrm{ss}, \delta \hat{e}(\theta^\mathrm{u}), g^\mathrm{ss}(q_+, \mu, c)).
\end{align}

Finally we observe that the location of the fronts and the symmetry operators imply properties of the function $g$. By assumption the location of fronts are zeros of the component $g^\mathrm{u}$ as we require $v^{u} = 0$ and $v^\mathrm{uu} = 0$ to lie in the center-stable manifold of the Turing family. Thus, as a consequence of Hypothesis~\ref{hyp:global} 
\changed{\begin{align}
    g^\mathrm{u} (s, \kappa^*(s), \alpha, 0, z(s), 0) = 0, \label{eqn:g-zeros}
\end{align}}
for all \changed{$\alpha,s \in S^1$}. Similarly by acting on the point $p$ in \eqref{eqn:pW-} by the equivariance $\Gamma_\alpha$ we obtain
\changed{\begin{align}
    g^\mathrm{ss}(s, \kappa^*(s), \alpha, 0, z(s), c) &= T_\alpha v^\mathrm{ss}_*(s), \label{eqn:g-equivariance}
\end{align}}
for all \changed{$\alpha,s \in S^1$}. This observation will allow us to find a unique solution for fixed \changed{phase $\theta^\mathrm{s}$} and obtain the circle $S^1$ of its time translations by applying $\Gamma_\alpha$.

Finally we need to make non-degeneracy conditions on the bifurcation diagram of the fronts.
\begin{Hypothesis}[Non-Degeneracy Condition] \label{hyp:degeneracy}
    The function \changed{$\kappa^*(s)$} is bounded away from zero, and the Jacobian matrix \changed{$\rmD_{\mu,\kappa}g^\mathrm{u}(s, \alpha, \kappa^*(s), 0, z(s), 0)$} is invertible for all \changed{$\alpha, s \in  S^1$}. Also the skeleton function $z$ has the property that its global minima and maxima are unique.
\end{Hypothesis}

One additional assumption on the skeleton $z$ is required to construct asymmetric defect solutions. Asymmetric defects connect fronts from different locations on the skeleton and so we require that the skeleton $z$ have a property that allows this. We also require one additional invertibility condition to solve the matching equations.

\begin{Hypothesis}[Asymmetric Contact Defects] \label{hyp:asymmetric}
\changed{We assume that $s_1,s_2\in S^1$ satisfy $z(s_1)=z(s_2)$ and $z'(s_1)z'(s_2)\neq0$ (note that $s_1,s_2$ are then locally part of a smooth one-parameter family parametrized by $(s_2=\mathcal{Z}(s_1)$ for some smooth function $\mathcal{Z}$). We also assume that the derivative $\rmD_{\varphi,\kappa}g^\mathrm{u}(s, \alpha, \kappa^*(s), 0, z(s), 0)$ is invertible for all $\alpha,s\in S^1$.}
\end{Hypothesis}

Note this hypothesis is effectively a non-degeneracy condition which is satisfied away from local maxima and minima of $z$. We can now state the fully rigorous results.

\begin{Theorem}[Symmetric Contact Defects] \label{thm:symmetric_contact}
Suppose Hypotheses~\ref{hyp:resolvent},~\ref{hyp:background},~\ref{hyp:Turing},~\ref{hyp:global}, and~\ref{hyp:degeneracy} hold. In particular, there are periodic functions $z:S^1 \to J$ and $\kappa^*: S^1 \rightarrow K$ for interval $J\subset \R$ and $K \subset \R^+$ describing the bifurcation diagram of the front and the core Turing wavenumber selected by each front respectively. Then for every $n\in\Z$ large enough, \changed{$s\in S^1$}, and $\varphi_0 \in \{ 0, \pi\}$ there exists a unique symmetric defect solution with temporal frequency $\omega = \omega_0(\mu)$ and core region of width $L=(2\pi n + \varphi_0-\changed{s})/\kappa^*(\changed{s})$ for parameter value $\mu = z(\changed{s}) + \bar{\mu}$ where $\bar{\mu} = \mathcal{O}(\mathrm{e}^{-\eta L})$ for some $\eta > 0$. 
\end{Theorem}

Here we do not refer to spirals and targets explicitly as they are solutions with distinct reversible symmetries (see Section~\ref{section:set_up_and_hyp}). Instead, Theorem~\ref{thm:symmetric_contact} can be applied to both reversers separately, as discussed in Section~\ref{section:proof_symmetric_snaking}, to conclude the presence of both spirals and targets. We also note that the contact defects found will have the same temporal frequency as the underlying front solutions which is equal to the temporal frequency of the homogeneous oscillations. Otherwise the theorem is essentially the same, but applies to any reaction-diffusion equation satisfying the listed hypotheses. 

\begin{Theorem}[Asymmetric Contact Defects] \label{thm:asymmetric_contact}
Suppose Hypotheses~\ref{hyp:resolvent},~\ref{hyp:background},~\ref{hyp:Turing},~\ref{hyp:global},~\ref{hyp:degeneracy}, and~\ref{hyp:asymmetric} hold. \changed{For every $n\in \N$ large enough, $s_1,s_2$ as in Hypothesis~\ref{hyp:asymmetric}, and $\theta\in  S^1$, there exists a unique asymmetric defect solution at parameter value $\mu=z(s_1)+\bar\mu$ with $\bar\mu=\bigo(|\Delta\kappa|/L+\mathrm{e}^{-\eta L})$ whose core region has width $L$ given by \eqref{eqn:asym_core_width} and whose asymptotic wave train profiles have a phase difference of $\theta$. These defects will travel with speed $c_\mathrm{d}(\Delta\kappa, \alpha) = \Delta \kappa / 2L  + \mathcal{O}(|\Delta\kappa|/L^2+\mathrm{e}^{- \eta L})$ and have temporal frequency $\omega = \omega_0(\mu) + \bigo(c_\mathrm{d}^2)$ where $\Delta \kappa = \kappa(s_2) - \kappa(s_1)$ is the difference of the wavenumbers of the Turing patterns selected by the front and back, and $\omega_0(\mu)$ is the temporal frequency of the homogeneous oscillation. The background oscillations consist of wave trains with spatial wavenumber $k = \bigo(c_\mathrm{d})$.}
\end{Theorem}

\section{Proof of Theorem~\ref{thm:symmetric_contact}}
\label{section:proof_symmetric_snaking}

Our strategy for this proof is to use the periodicity of the skeleton function $z$ to find symmetric defects parameterized by \changed{$s\in  S^1$} corresponding to their phase coordinate at the midpoint in the neighborhood $\mathcal{V}$. The major drawback of this approach is that we do not directly show that a contact defect exists for each value of $\mu$. Instead this is implicit in our result through our parameterization of the bifurcation diagram $z(\varphi)$. For any value of \changed{$\mu = z(s)$} with a front, we use the construction below to find nearby defect solutions for $\mu = z(\changed{s}) + \bar\mu$ where $\bar\mu$ is shown to be exponentially small. Then, so long as we are away from the saddle-nodes of the fronts, and using the assumption that $\rmD_{\mu,\kappa}g^\mathrm{u}$ is invertible, we can use the Implicit Function Theorem to write the bifurcation curve as a graph over $\mu$ close to $z(\changed{s})$. By uniqueness of these solutions, we can patch these graphs for different values of \changed{$s$} together to reassemble the symmetric branches of the diagram as functions of $\mu$.

The advantage of this method is that we naturally obtain the leading order periodic structure of the snaking diagram from the skeleton function $z(\changed{s})$, including near its global extrema. 

Our objective is to find an exchange solution $v(x)$ through the neighborhood $\mathcal{V}$ that connects the center-unstable and center-stable manifolds of the background states which is also $\mathcal{R}$-symmetric. Using this symmetry condition, it suffices for us to solve the following problem: find an exchange solution $v(x)$ such that,
\begin{align}
\begin{split}
	&v(x) \in \mathcal{V} \: \text{for all} \: x\in[-L,L], \\
	&v(-L) \in \rmW^\mathrm{cu}(\changed{\mathcal{E}(\mu)})\cap \Sigma^\mathrm{s}, \\
	&v(0) \in \text{Fix} \mathcal{R}.
 \end{split}
\end{align}

From Lemma~\ref{lemma:fenichel} we know the action of the reverser $\mathcal{R}$ on our Fenichel coordinates and have exact criteria for when the exchange solution $v(x) = v(x; \mu, \kappa_0, \varphi_0, \theta^\mathrm{s}_0, \theta^\mathrm{u}_0, v^\mathrm{ss}_0, v^\mathrm{uu}_0)$ is $\mathcal{R}$-symmetric. This is the case if, and only if, $\varphi_0 \in \{0,\pi\}$, $\theta^\mathrm{s}_0 = \theta^\mathrm{u}_0$ (mod $2\pi$), $v_0^\mathrm{ss} = v_0^\mathrm{uu}$ and $c=0$. So for the remainder of this section we restrict to those exchange solutions $v(x;\mu, 0, \kappa_0, \varphi_0, \theta_0, \theta_0, v_0, v_0)$ where $\varphi_0 \in \{0, 2\pi\}$. Since $c$ is forced to be zero, $\omega = \omega_0$ and we ignore the dependence on the \changed{wave speed}.

This reduces the problem to solving the matching condition $v(-L) \in \rmW^\mathrm{cu}(\changed{\mathcal{E}(\mu)})\cap\Sigma^\mathrm{s}$. Lemma~\ref{lemma:exchange} shows that $v(-L)$ can be written as
\begin{align*}
    & v(-L; \mu, \kappa_0, \varphi_0, \theta_0, \theta_0, v_0, v_0) \\
    & = (\varphi_0-\kappa_0 L +\expsmall, \kappa_0+\expsmall, \delta \hat{e}(\theta_0), v_0 + \expsmall, \expsmall, \expsmall).
\end{align*}
Hypothesis~\ref{hyp:global} implies that elements of manifold $\rmW^\mathrm{cu}(\changed{\mathcal{E}(\mu)})\cap\Sigma^\mathrm{s}$ are of the form $(\varphi, \kappa, \delta \mathrm{e}^{i\theta^\mathrm{s}}, g^\mathrm{ss}(q), g^\mathrm{u}(q), v^\mathrm{uu})$ for $q = (\varphi, \kappa, \theta^\mathrm{s}, v^\mathrm{uu}, \mu)$. Thus, solving the matching condition $v(-L) \in \rmW^\mathrm{cu}(\changed{\mathcal{E}(\mu)})\cap\Sigma^\mathrm{s}$ is equivalent to constructing solutions $(L,\varphi, \theta^\mathrm{s}, \kappa, v^\mathrm{uu}, \mu, \kappa_0, \theta_0, v_0)$  of the system
\begin{align}
    (\varphi_0-\kappa_0 L +\expsmall, \kappa_0+\expsmall, \delta \hat{e}(\theta_0), v_0 + \expsmall, &\expsmall, \expsmall) \nonumber \\
    = (\varphi, \kappa, \delta \hat{e}(\theta^\mathrm{s}), g^\mathrm{ss}(q), g^\mathrm{u}(q), v^\mathrm{uu})&.
\end{align}
Since the existence statement should be true for all $L$ sufficiently large, we will look for solutions in the limit as $L$ grows to infinity.

We immediately conclude that $\theta^\mathrm{s} = \theta_0 =: \changed{\alpha}$ is a free parameter unfolding the $S^1$-equivariant symmetry of our solutions. By fixing \changed{$\alpha$} we factor out this symmetry so that we are looking for unique solutions. We can also solve for $v^\mathrm{uu}$ using the Implicit Function Theorem in the limit as $L\rightarrow \infty$ and find that $v^\mathrm{uu} = \expsmall$. The remaining equations are:
\begin{align}
	&\varphi = \varphi_0 - \kappa_0 L + \expsmall  \quad \mod{2\pi}, \\
	&\kappa = \kappa_0 + \expsmall, \\
	&g^\mathrm{ss}(\varphi, \kappa, \alpha, \expsmall, \mu) = v_0 + \expsmall,\\
	&g^{u}(\varphi, \kappa, \alpha, \expsmall, \mu) = \expsmall,
\end{align}
for the unknowns $(\varphi, \kappa, \kappa_0, \mu, L, v_0)$. We are looking for a 1-parameter family of solutions for each $\alpha\in  S^1$ (we have factored out the equivariance), so this dimension count makes sense. We also have a choice of which variable to use to parameterize the family. The intuitive choice is $L$ as this is our proxy for the width of the core region, however it will prove to be easier if we instead use \changed{$\varphi = s$ as this is the branch parameter for the skeleton function $z$}. With this in mind we relabel our unknowns to make use of Hypothesis~\ref{hyp:degeneracy} and express everything as functions of \changed{$\varphi = s$. Let $\mu = z(s) + \bar\mu$, $\kappa = \kappa^*(s) + \bar\kappa$, $\kappa_0 = \kappa^*(s) + \bar\kappa_0$ and $v_0 = T_{\alpha} v^\mathrm{ss}_*(s) + \bar{v}$ so that},
\begin{align}
	&\changed{s} = \varphi_0 - (\kappa^*(\changed{s}) + \bar\kappa) L + \expsmall \quad \mod{2\pi},  \\
	& \bar\kappa = \bar\kappa_0 + \expsmall, \\
    &g^\mathrm{ss}(\changed{s}, \kappa^*(\changed{s}) + \bar\kappa, \alpha, \expsmall, z(\changed{s}) + \bar\mu) = T_\alpha v^\mathrm{ss}_*(\changed{s}) + \bar{v} +  \expsmall, \\
	&g^{u}(\changed{s}, \kappa^*(\varphi) + \bar\kappa, \alpha, \expsmall, z(\changed{s}) + \bar\mu) =  \expsmall.
\end{align}

Recall equations \eqref{eqn:g-zeros} and \eqref{eqn:g-equivariance},
\begin{align*}
    g^\mathrm{u} (\changed{s}, \kappa^*(\changed{s}), \alpha, 0, z(\changed{s}), 0) &= 0, \\
    g^\mathrm{ss}(\changed{s}, \kappa^*(\changed{s}), \alpha, 0, z(\changed{s}), 0) &= T_\alpha v^\mathrm{ss}_*(\changed{s}),
\end{align*}
\changed{for all $\alpha, s \in S^1$}, providing us known values of $g^\mathrm{u}$ and $g^\mathrm{ss}$ about which we can expand to find
\begin{align}
	&\changed{s} = \varphi_0 - (\kappa^*(\changed{s}) + \bar\kappa) L + \expsmall \quad \mod{2\pi} \label{eqn:sym_match_1_a}\\
	& \bar\kappa - \bar\kappa_0 =  \expsmall \label{eqn:sym_match_1_b} \\
    &\bar{v} + \rmD_{\mu,\kappa} g^\mathrm{ss}(\changed{s}, \kappa^*(\changed{s}), \alpha, 0, z(\changed{s})) \begin{pmatrix} \bar\mu \\ \bar\kappa \end{pmatrix} =  R_1(\bar\mu,\bar\kappa) + \expsmall \label{eqn:sym_match_1_c} \\
	&\rmD_{\mu,\kappa} g^{u}(\changed{s}, \kappa^*(\changed{s}), \alpha, 0, z(\changed{s})) \begin{pmatrix} \bar\mu \\ \bar\kappa \end{pmatrix} = R_2(\bar\mu, \bar\kappa) + \expsmall \label{eqn:sym_match_1_d}
\end{align}
where $R_1$ and $R_2$ are quadratic in $(\bar\mu, \bar\kappa)$ by Taylor's theorem. We have assumed that the derivatives of $g$ are uniformly bounded meaning the other first derivative terms contribute at most $\expsmall$ and can be absorbed into the higher order terms. We are now looking for a solution $(L, \bar\kappa, \bar{\kappa}_0, \bar{v}, \bar\mu)$ as functions of \changed{$s$}. 

Hypothesis~\ref{hyp:degeneracy} implies that $\rmD_{\mu,\kappa} g^{u}:\R^2\to\R^2$ is invertible (recall that the weak unstable subspace is two-dimensional so $g^\mathrm{u}$ takes values in $\R^2$). From this we conclude that the linear operator of $(\bar\mu,\bar\kappa,\bar\kappa_0,\bar{v})$ on the left-hand side of the final three equations, \eqref{eqn:sym_match_1_b}-\eqref{eqn:sym_match_1_d}, is invertible, and we use the Implicit Function Theorem to solve uniquely for $(\bar\mu,\bar\kappa,\bar\kappa_0,\bar{v})$ along with the estimates $(\bar\mu,\bar\kappa,\bar\kappa_0,\bar{v})=\expsmall$. All that remains is to solve the first equation, which, upon substituting $\bar\kappa=\expsmall$, is given by
\begin{align}
    \changed{s} = \varphi_0 - \kappa^*(\changed{s}) L +  \expsmall \quad \mod{2\pi}.
\end{align}
Since this is an equation posed on $\R/2\pi\Z$, we write $L = (2\pi n + l)/\kappa^*(\changed{s})$ for $n\in \mathbb{Z}$ and $l\in S^1$ where we note that we assumed in Hypothesis~\ref{hyp:degeneracy} that $\kappa^*(\changed{s})\neq0$. The resulting equation for $l$ is given by
\begin{align}
    \changed{s} = \varphi_0 - l + \mathcal{O}(\mathrm{e}^{-\bar\eta n}),
\end{align}
which we can solve for $l = \varphi_0 - \changed{s} + \mathcal{O}(\mathrm{e}^{-\bar\eta n})$ as a function of \changed{$s$} and $n$ (recall that $\varphi_0\in\{0,\pi\}$ is fixed) so that
\begin{equation}
    L = \frac{\varphi_0 - \changed{s} + 2\pi n}{\kappa^*(\changed{s})} + \mathcal{O}(\mathrm{e}^{-\bar\eta n})
\end{equation}
for each $n \in \mathbb{Z}$ large enough. That is, for each front in the bifurcation diagram $\mathcal{B}$ there exist a countably infinite number of defect solutions nearby, indexed by a natural number $n$, and which spend an increasing amount of time near the Turing family. This proves Theorem~\ref{thm:symmetric_contact}.

In this proof we reduced the problem to a single matching condition using criteria for $\mathcal{R}$-symmetric solutions in $\mathcal{V}$. We can equally determine such a reduction for $\Gamma_\pi\mathcal{R}$-symmetric solutions. Indeed, an exchange solution $v$ is $\Gamma_\pi\mathcal{R}$-symmetric if, and only if, $\varphi_0 \in \{0,\pi\}$, $\theta^\mathrm{s}_0 - \theta^\mathrm{u}_0 = \pi \mod{2\pi}$, $v^\mathrm{ss}_0 = T_\pi v^\mathrm{uu}_0$ and $c=0$, using the Exchange Lemma asymptotics. Once this identification is made, the remainder of the proof continues in exactly the same way and we find that we have branches of both target and spiral defects as claimed in Section~\ref{section:set_up_and_hyp}.

\section{Proof of Theorem~\ref{thm:asymmetric_contact}}
\label{section:asymmetric_snaking}

This proof follows in much the same way as the previous, except that we cannot use symmetry to reduce the matching equations. As a consequence, there are additional degrees of freedom that we need to take account of. 

Most importantly, we will need to consider the equations in a comoving frame with a non-zero wave speed $c$. As discussed in Section~\ref{section:spatial_dynamics}, we need to be careful in our construction that we end up with contact defects at the end. Recall from Lemma~\ref{lemma:group_velocities}, for all $c$ small enough we can find wave trains with spatial wavenumber $k = k(c)$ and temporal frequency $\omega = \omega(c)$ which have group velocity $c$. \changed{As noted in Section~\ref{section:set_up_and_hyp}, by} setting $\omega = \omega(c)$ we ensure that any defect we find will be a contact defect and we drop the dependence of our manifolds on $\omega$.

We will also find solutions that are parameterized both by \changed{the branch parameter $\varphi = s$ and a phase offset $\theta = \theta^\mathrm{u}-\theta^\mathrm{s}$}. The \changed{branch parameter $s$} will trace out the winding of the snaking diagram as it did in the symmetric case, \changed{while $\theta$} will parameterize the circular cross-sections of the tubular neighborhoods. Note this is distinct from the symmetric case where $\theta^\mathrm{s} = \changed{\alpha}$ parameterized the $S^1$-equivariance of our solutions for fixed $\mu$. In this case we can again set $\theta^\mathrm{s} = \alpha$, which parameterizes this equivariance, but we can also still vary $\theta^\mathrm{u}$ to make all possible background phase offsets. \changed{In the proof we do not solve for either angle $\alpha$ or $\theta$, instead we show that for any fixed choices of these angles, we can find a unique contact defect.}

Again we wish to find an exchange solution $v(x)$ through $\mathcal{V}$ that connects the center-unstable and center stable manifolds of the background states. However, this time we cannot simplify the system using reversibility. Instead, we are looking for a solution $(L,c,\varphi, \varphi_-, \varphi_+,\kappa, \kappa_-, \kappa_+, \theta^\mathrm{s}, \theta_0^\mathrm{s}, \theta^\mathrm{u}, \theta_0^\mathrm{u}, v_0^\mathrm{ss}, v_+^\mathrm{ss}, v_0^\mathrm{uu}, v_-^\mathrm{uu})$    to the two matching conditions
\begin{align}
    (\varphi, \kappa, \delta \hat{e}(\theta_0^\mathrm{s}), v_0^\mathrm{ss}, \expsmall, \expsmall) = (\varphi_-, \kappa_-, \delta \hat{e}(\theta^\mathrm{s}), g^\mathrm{ss}(q_-), g^\mathrm{u}(q_-), v_-^\mathrm{uu}), \\
    (\varphi + 2\kappa L + 2cL^2 + \expsmall, \kappa + 2cL + \expsmall, \expsmall, \expsmall, \delta \hat{e}(\theta_0^\mathrm{u}), v_0^\mathrm{uu})& \nonumber \\
    = (-\varphi_+, \kappa_+, g^\mathrm{u}(q_+), v_+^\mathrm{ss}, &\delta \hat{e}(\theta^\mathrm{u}), g^\mathrm{ss}(q_+)),
\end{align}
where $q_- = (\varphi_-, \kappa_-, \theta^\mathrm{s}, v_-^\mathrm{uu}, \mu, c)$ and $q_+ = (\varphi_+,  \kappa_+, \theta^\mathrm{u}, v_+^\mathrm{ss}, \mu, c)$. The exchange asymptotics have been written slightly differently here than was presented in Lemma~\ref{lemma:exchange}. The leading order behavior of the center coordinates, $\varphi$ and $\kappa$, have been shifted so that the initial value is found at $x=-L$ instead of $x=0$. This is an aesthetic choice which makes the matching conditions at $x=-L$ simpler, so that we can concentrate on solving the matching conditions at $x=L$.

We immediately have that $\varphi = \varphi_-$, $\kappa = \kappa_-$, $\theta_0^\mathrm{s} = \theta^\mathrm{s}$, $\theta_0^\mathrm{u} = \theta^\mathrm{u}$ and we can use the Implicit Function Theorem to solve for $v_-^\mathrm{uu} = \expsmall$ and $v_+^\mathrm{ss} = \expsmall$. The remaining equations are then
\begin{align}
    -\varphi_+ = \varphi_- + 2cL^2 + 2\kappa_- L \ + &\ \expsmall \quad \mod{2\pi}, \\
    \kappa_+ = \kappa_- + 2cL \ + &\ \expsmall, \\
    g^\mathrm{u}(\varphi_-, \kappa_-, \theta^\mathrm{s}, \expsmall, \mu, c) &= \expsmall, \\
    g^\mathrm{u}(\varphi_+, \kappa_+, \theta^\mathrm{u}, \expsmall, \mu, c) &= \expsmall, \\
    g^\mathrm{ss}(\varphi_-, \kappa_-, \theta^\mathrm{s}, \expsmall, \mu, c) &= v_0^\mathrm{ss}, \\
    g^\mathrm{ss}(\varphi_+, \kappa_+, \theta^\mathrm{u}, \expsmall, \mu, c) &= v_0^\mathrm{uu}.
\end{align}

Excluding the final two equations which are infinite dimensional, we have a 6-dimensional set of equations in the seven unknowns $\varphi_\pm$, $\kappa_\pm$, $c$, $L$ and $\mu$, meaning we will have to solve for at least one of $\varphi_\pm$. \changed{Using Hypothesis~\ref{hyp:asymmetric} we let $\varphi_- = s$ to describe the branch parameter and let $\varphi_+ = \mathcal{Z}(s) + \bar{\varphi}_+$ with the aim of solving for the correction $\bar{\varphi}_+$}. Suppose we choose \changed{$\varphi_- = s$} away from the global maxima and minima of $z$. Then like the symmetric case we look to show the resulting solutions are formed at a small perturbation of the periodic skeleton determined by \changed{$\varphi_- = s$ and $\varphi_+ = \mathcal{Z}(s)$. Let $\kappa_- = \kappa^*(s) + \bar{\kappa}_-$, $\kappa_+ = \kappa^*(\mathcal{Z}(s)) + \bar{\kappa}_+$, $\mu = z(s) + \bar{\mu} = z(\mathcal{Z}(s)) + \bar{\mu}$, $v_0^\mathrm{ss}  = T_{\theta^\mathrm{s}} (v_*^\mathrm{ss}(s) + \bar{v}_0^\mathrm{ss})$ and $v_0^\mathrm{uu}  = T_{\theta^\mathrm{u}} (v_*^\mathrm{ss}(\mathcal{Z}(s)) + \bar{v}_0^\mathrm{uu})$.}  Substituting these we get
\begin{align}
    -\mathcal{Z}(\changed{s}) - \bar{\varphi}_+ = \changed{s} + 2cL^2 + 2(\kappa^*(\changed{s}) + \bar{\kappa}_-)L \ + \expsmall \mod{2\pi},& \label{eqn:asymm_match_L} \\
    \kappa^*(\mathcal{Z}(\changed{s})) + \bar{\kappa}_+ = \kappa^*(\changed{s}) + \bar{\kappa}_- + 2cL \ + \expsmall, & \\
    g^\mathrm{u}(\changed{s}, \kappa^*(\changed{s}) + \bar{\kappa}_-, \theta^\mathrm{s}, \expsmall, z(\changed{s}) + \bar{\mu}, c) = \expsmall, & \\
    g^\mathrm{u}(\mathcal{Z}(\changed{s}) + \bar{\varphi}_+, \kappa^*(\mathcal{Z}(\changed{s})) + \bar{\kappa}_+, \theta^\mathrm{u}, \expsmall, z(\mathcal{Z}(\changed{s})) + \bar{\mu}, c) = \expsmall, & \\
    g^\mathrm{ss}(\changed{s}, \kappa^*(\changed{s}) + \bar{\kappa}_-, \theta^\mathrm{s}, \expsmall, z(\changed{s}) + \bar{\mu}, c),
    = T_{\theta^\mathrm{s}} (v_*^\mathrm{ss}(\changed{s}) + \bar{v}_0^\mathrm{ss}), & \label{eqn:asymm_match_vss} \\
    g^\mathrm{ss}(\mathcal{Z}(\changed{s}) + \bar{\varphi}_+, \kappa^*(\mathcal{Z}(\changed{s})) + \bar{\kappa}_+, \theta^\mathrm{u}, \expsmall, z(\mathcal{Z}(\changed{s})) + \bar{\mu}, c), 
    =T_{\theta^\mathrm{u}} (v_*^\mathrm{ss}(\mathcal{Z}(\changed{s})) + \bar{v}_0^\mathrm{uu}).& \label{eqn:asymm_match_vuu}
\end{align}

We start by solving the second, third and fourth equations. By expanding the third equation about \changed{$q_- = (s, \kappa^*(s), \theta^{s}, 0, z(s),0)$} and the fourth equation about \changed{$q_+ = (\mathcal{Z}(s),  \kappa^*(\mathcal{Z}(s)), \theta^\mathrm{u}, 0, z(\mathcal{Z}(s)), 0)$} and rearranging to put the leading order linear operators on the left-hand side, these equations become
\begin{align}
    \rmD_{\varphi,\kappa}g^\mathrm{u}(q_+)\begin{pmatrix} \bar{\varphi}_+ \\ \bar{\kappa}_+ \end{pmatrix} + \frac{\partial g^\mathrm{u}}{\partial \mu}(q_+) \bar{\mu} + \frac{\partial g^\mathrm{u}}{\partial c}(q_+) c + R_2(\bar{\varphi}_+,\bar{\kappa}_+, \bar{\mu}, c) + \expsmall = 0,& \label{eqn:asymm_matching_expanded_st}\\
    \rmD_{\kappa,\mu} g^\mathrm{u} (q_-) \begin{pmatrix} \bar{\kappa}_- \\ \bar\mu \end{pmatrix} + \frac{\partial g^\mathrm{u}}{\partial c}(q_-) c + R_1(\bar{\kappa}_-, \bar{\mu}, c) +  \expsmall = 0,& \\
    \bar{\kappa}_+ - \bar{\kappa}_- - 2Lc + \expsmall = \kappa^*(\changed{s}) - \kappa^*(\mathcal{Z}(\changed{s})),& \label{eqn:asymm_matching_expanded_end}
\end{align}
where the $\varphi, \kappa$ and $\mu$ subscripts denote partial derivatives, and the functions $R_1$ and $R_2$ are quadratic in their arguments by Taylor's theorem. Again we use that the derivatives of $g^\mathrm{u}$ are uniformly bounded in Hypothesis~\ref{hyp:degeneracy} to subsume other derivative terms into the order $\expsmall$ terms. Our goal is to solve these equations for the unknowns $(\bar\varphi_+,\bar\kappa_+, \bar\kappa_-,\bar\mu,c)$ by the Implicit Function Theorem which requires that the leading order linear operator on the left-hand side be invertible. We can write this operator as
\begin{align}
    \left( \begin{array}{@{}cc|c@{}}
        \rmD_{\varphi,\kappa}g^\mathrm{u}(q_+) & \begin{matrix}
            \begin{matrix}
                0 \\ 0
            \end{matrix} & \frac{\partial g^\mathrm{u}}{\partial \mu}(q_+)
        \end{matrix} & \frac{\partial g^\mathrm{u}}{\partial c}(q_+) \\
        \bigzero & \rmD_{\kappa,\mu} g^\mathrm{u} (q_-) & \frac{\partial g^\mathrm{u}}{\partial c}(q_-) \\ 
        \hline
        \begin{matrix}
            0\hphantom{00} & 1
        \end{matrix} & \begin{matrix}
            -1\hphantom{0000} & 0
        \end{matrix} & -2L
    \end{array} \right) =: \left( \begin{array}{c|c}
         A & b \\
         \hline
         \begin{matrix}
             0 & 1 & -1 & 0
         \end{matrix} & -2L 
    \end{array}\right) \label{eqn:jacobian_for_solve}
\end{align}

The main $4\times 4$ block, $A$, is invertible by Hypotheses~\ref{hyp:degeneracy} and~\ref{hyp:asymmetric} where we assume that its diagonal blocks are invertible. To show the full matrix is invertible, we show that it has a trivial kernel for all $L$ sufficiently large. Suppose $(\bar\varphi_+,\bar\kappa_+, \bar\kappa_-,\bar\mu,c)$ is in the kernel of the operator, then by matrix multiplication and using the invertibility of $A$ we find
\begin{align}
    \begin{bmatrix}
        \bar\varphi_+ \\ \bar\kappa_+ \\ \bar\kappa_- \\ \bar\mu
    \end{bmatrix} = -c A^{-1} b, \quad \bar\kappa_+ - \bar\kappa_- - 2Lc = 0.
\end{align}
Using the first equation we can eliminate $\bar\kappa_\pm$ from the second equation which yields
\begin{equation}
    c ( \langle A^{-1}_3 - A^{-1}_2, b \rangle - 2L ) = 0
\end{equation}
where $A^{-1}_j$ is the $j$-th row of the matrix $A^{-1}$ and $\langle \cdot, \cdot \rangle$ is the vector dot product. Since, $A$ and $b$ are made up of derivatives of $g^\mathrm{u}$, which we assume are uniformly bounded for all sufficiently large $L$, the only solution to this equation is $c = 0$. This implies the vector $(\bar\varphi_+,\bar\kappa_+, \bar\kappa_-,\bar\mu,c)$ is exactly zero.

Applying the Implicit Function Theorem, we can find unique solutions $(\bar\varphi_+,\bar\kappa_+, \bar\kappa_-,\bar\mu,c)$ to equations \eqref{eqn:asymm_matching_expanded_st}-\eqref{eqn:asymm_matching_expanded_end} which have expansions
\begin{align}
    c &= \frac{\Delta\kappa}{2L} + \changed{\bigo\left(\frac{|\Delta\kappa|}{L^2} + \mathrm{e}^{-\eta L}\right)}, \\
    (\bar\varphi_+,\bar\kappa_+, \bar\kappa_-,\bar\mu) &= \bigo\left(\frac{|\Delta\kappa|}{L} + \mathrm{e}^{-\eta L}\right),
\end{align}
where $\Delta\kappa := \kappa^*(\mathcal{Z}(\changed{s})) - \kappa^*(\changed{s})$ is the difference between the selected wavenumbers. Using these expansions we can solve equations \eqref{eqn:asymm_match_vss} and \eqref{eqn:asymm_match_vuu} and again find that $\bar{v}_0^\mathrm{ss}, \bar{v}_0^\mathrm{uu} = \bigo(|\Delta\kappa|/L + \mathrm{e}^{-\eta L})$.

All that remains is to solve the equation \eqref{eqn:asymm_match_L} for $L$. Using our new asymptotics for $\bar\kappa_+$ and $\bar\kappa_-$ this equation reduces to 
\begin{align}
    -\mathcal{Z}(\changed{s}) = \changed{s} + B\Delta\kappa + L(\kappa^*(\mathcal{Z}(\changed{s})) + \kappa^*(\changed{s})) + \bigo\left(\frac{|\Delta\kappa|}{L} + \mathrm{e}^{-\eta L}\right) \mod{2\pi}
\end{align}
for some uniformly bounded constant $B$ which is determined by the entries of the inverse of the matrix defined in \eqref{eqn:jacobian_for_solve}. Let $L = (l + 2\pi n)/(\kappa^*(\mathcal{Z}(s)) + \kappa^*(s))$ for $n\in\Z$ large and $l\in  S^1$, then by the Implicit Function Theorem we conclude
\begin{align}
    L = \frac{-(\mathcal{Z}(\changed{s}) + \changed{s} + B\Delta\kappa) + 2\pi n}{\kappa^*(\mathcal{Z}(\changed{s})) + \kappa^*(\changed{s}) } + \mathcal{O}\left(\frac{|\Delta\kappa|}{n} + \mathrm{e}^{-\eta n}\right)
\end{align}
completing the proof.

This proof only depends on the choice front and back selected assuming we can find valid choices for $\varphi_\pm$ \changed{through $s$ and $\mathcal{Z}(s)$}. We can exactly recover our result for symmetric contact defects by choosing $\varphi_- = \varphi_+$ for which $\Delta\kappa = 0$, which reproduces both snaking curves corresponding to $\varphi_0 = 0,\pi$. Along with our symmetric contact defects, for each value of $\mu \in J$ we get an $ S^1$-family of contact defects with a phase offset \changed{$\theta\in  S^1$} in their background oscillations which move at very small wave speeds of the order $\bigo(\mathrm{e}^{-\eta L})$. 

There are also fully asymmetric solutions formed by taking distinct values $\varphi_-$ and $\varphi_+$ which select different wavenumbers. In this case, the defects move at larger speeds $c \approx \Delta\kappa/2L$. They also have asymmetric core regions which have core Turing regions shifted in $x$ and modulated in spatial wavenumber $\kappa$. These branches also differ from the regular snaking situation, as they will not lie exponentially close to the position predicted by the skeleton $z$. Instead, they are a larger $\Delta\kappa/L$ perturbation, which is small for sufficiently large $L$ and for $\mu$ close to the symmetric branch saddle-nodes. 

\section{Discussion} \label{section:discussion}

In this paper, we have built a rigorous framework to understand the homoclinic snaking of contact defects in reaction-diffusion equations. We have shown that the bifurcation diagram of contact defects can be predicted qualitatively and quantitatively by the bifurcation diagram of their generating heteroclinic front. We found that these diagrams are very different from those found in snaking bifurcations of stationary localized patterns due to an extra degree of freedom inherited from an additional center-direction in the linearization of their spatially homogeneous background oscillations. As a consequence, we demonstrate the existence of contact-defect solutions with arbitrary background phase offsets, a fact that has not been observed previously in the literature. Moreover, we argue that this is an essential feature of the snaking of defect solutions as it allows us to reconcile the robustness properties of contact defects with the theory of homoclinic and heteroclinic bifurcations.

\changed{We note that, while we expect that our nondegeneracy conditions hold for typical systems, we cannot prove that the Brussellator system (or any other pattern-forming reaction-diffusion system) satisfies the hypotheses we need for our construction of contact defects. This is a common problem, and we mention the lack of a rigorous verification of the hypotheses stated in \cite{beck2009} for localized patterns in the Swift--Hohenberg equation as an additional example.}

We now mention a number of avenues for further study of homoclinic snaking in reaction-diffusion equations that follow naturally from the analysis presented here. 

\begin{figure}
    \centering
    \includegraphics[width=0.7\linewidth]{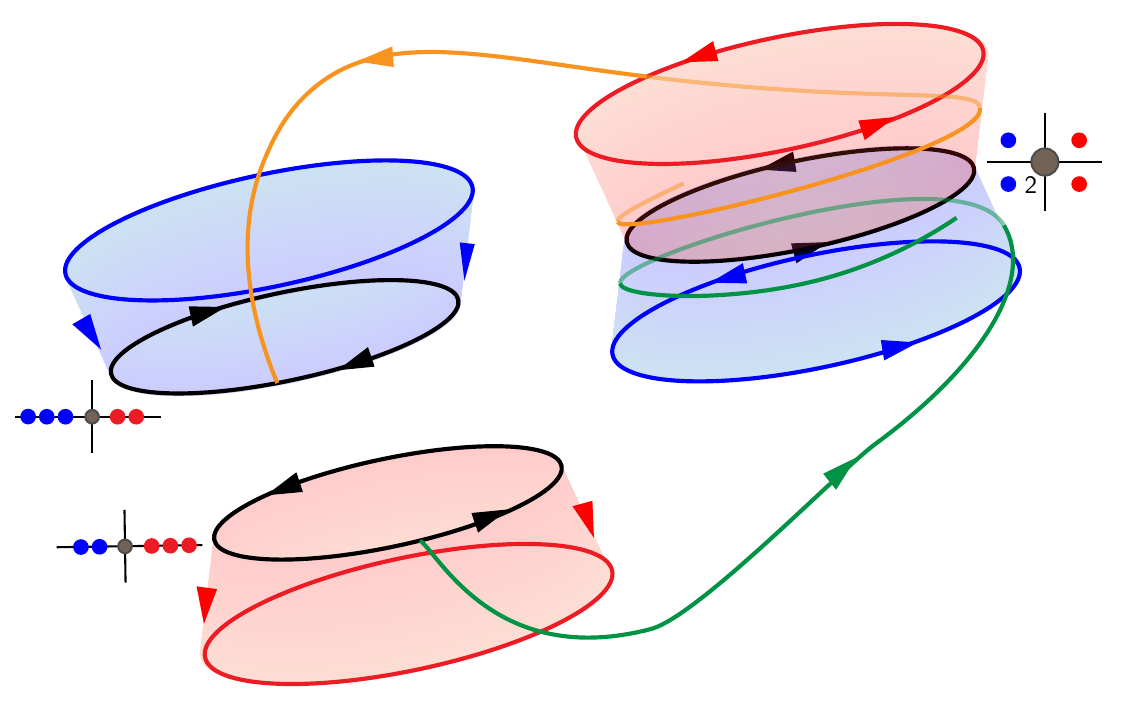}
    \caption{Schematic diagram of phase space for Source Defects. Background states are wave trains with nonzero group velocity relative to the defect core and have less center-stable directions than in the contact defect case.}
    \label{fig:source_geometry}
\end{figure}

\paragraph{Homoclinic snaking for source defects.}As alluded to above, other types of defects can be found in reaction-diffusion equations \cite{sandstede2004}. Of particular interest to us are source defects that are, by definition, asymptotic to wave trains with outwardly directed group velocities. Source defects differ in several ways from the contact defects we considered in this paper. The background states are now a pair of true wave trains with a non-zero group velocity directed outward relative to the defect core. A schematic of the resulting phase space geometry is pictured in Figure~\ref{fig:source_geometry}. A direct result of Lemma~\ref{lemma:spectra_wave trains} and the spatial Floquet-multipliers pictured in Figure~\ref{fig:spatial_spectrum_wave trains} is that the center-unstable manifold of the wave train at $x=-\infty$ has one dimension less than in the contact-defect case. Consequently, we should no longer expect that the heteroclinic front solutions will be formed generically. This can be rectified by including an additional parameter in the problem, namely the temporal frequency $\omega$. We believe that, so long as one can understand the impact of this change, source defects will also be found in snaking diagrams. However, we also expect that the nature of these snaking diagrams will be different. Unlike contact defects, source defects will exist generically for discrete values of the temporal frequency $\omega = \omega_\mathrm{nl}(k)$, which will form 1-dimensional branches when continued in a bifurcation parameter $\mu$. As such, we do not expect to find source defects with arbitrary phase offsets. Instead, we expect that the background oscillations again couple very weakly \cite{chapman2009, kozyreff2006} allowing only target and spiral defects. This also means that the expected snaking bifurcations will look similar to those found in the case of stationary localized solutions in the Swift--Hohenberg equation. In future work, we plan to extend our current framework to give a full account snaking for source defects. 

\paragraph{Snaking of contact defects in the presence of secondary bifurcations.} Recent numerical studies by Al Saadi et al.\ \cite{alsaadi2024} have found evidence of homoclinic snaking for defect solutions in the Leslie--Gower and Gilad--Meron models. In comparison to the studies by Tzou et al.\ \cite{tzou2013}, they consider parameter regions where homoclinic snaking of localized solutions occurs simultaneously with the destabilization of the background state at a uniform Hopf bifurcation. They demonstrate that this instability is inherited by the localized solution snaking diagrams which produces branches of contact defects. Moreover, through a codimension-2 process and by varying a second bifurcation parameter, these Hopf bifurcations themselves collide and annihilate to form a full snaking diagram of contact defects. We note that in building our theory, we have deliberately avoided scenarios where secondary bifurcations occur. In particular, we have confined ourselves to the situation where the background oscillations are purely stable. It would be interesting to use the framework built here to analyze these more complex scenarios and determine what relationship the localized and contact defect snaking diagrams have close to these instabilities.

\paragraph{Numerical computation of contact defects.} In this work we have not presented any numerical computations for contact defects and their snaking diagrams, save for a replication of the original numerics performed by Tzou et al. \cite{tzou2013}. While it is relatively easy to compute symmetric contact defects, constructing a numerical scheme that can compute and continue asymmetric contact defects presents a major challenge. In order to construct such a scheme, and to be able to assess the accuracy of our theoretical predictions, we need to be able to compute defects on large spatial domains and enforce the correct convergence for large $x$. 

A natural choice for such a problem is a core-far field  decomposition as employed by Morrissey and Scheel \cite{morrissey2015}, and Lloyd and Scheel \cite{lloyd2017} for example. Such a scheme seeks to decompose the desired solution $u_\mathrm{d}$ as a sum $\chi^-  u_\mathrm{wt}^- + \chi^+  u_\mathrm{wt}^+ + w$ where $\chi^\pm$ are smooth cut-off functions with support for $x<-d$ and $x>d$ respectively, and $w$ is a localized interface connecting the background wave trains. The advantages of a core-far field decomposition is that it scales well in the large domain limit and also gives you direct control over the background asymptotics so long as you are able to correctly enforce localization of the interface function $w$.

Unfortunately the nature of contact defects complicates this approach in practice. We know that convergence to the background wave trains is algebraic and not exponential, as is presupposed in previous implementations of this scheme. We believe it is possible to overcome this problem by considering the core-far field decomposition on the correct Kondratiev space, as has been used in previous theoretical studies of striped patterns \cite{jaramillo2015}. However, Fredholm index calculations show that this has direct implications on which boundary conditions can be used in the scheme. We also know that there is no asymptotic phase for contact defects and instead the asymptotic far-field state is, to leading order, given by $ u_\mathrm{wt}(kx-\tau + \theta(x))$ where $\theta(x)$ is, to leading order, logarithmic for large $x$. This logarithmic phase correction cannot be ignored in the core-far field decomposition as without it there is no hope for localization of the interface solution $w$ on large domains. Accounting for this requires a very careful set up of the problem, which we will present in a follow-up paper to this theoretical investigation.

\paragraph{Linear and nonlinear stability of snaking contact defects.} Finally we note that the results presented here make no statement about the stability of the contact defects. Similar to the existence problem, the natural structure of contact defects as glued fronts and backs may make it amenable to study using techniques from dynamical systems similar to those used in  \cite{makrides2019}. However, the double spatial eigenvalue structure of the asymptotic homogeneous state and the resulting algebraic decay of the defect tails imply that the Evans function cannot be extended analytically across $\lambda=0$ \cite{sandstede2004a}. Thus, studying spectral stability will be challenging, and establishing linear stability using pointwise estimates will be even more complicated \cite{howard2006}.

\section{Acknowledgements}
The authors gratefully acknowledge helpful discussions with Arnd Scheel. Roberts was supported by the NSF under grant DMS-2106566. Sandstede was partially supported by the NSF under grants DMS-2038039 and DMS-2106566.


\bibliographystyle{abbrv}
\bibliography{Snaking_Defects}

\end{document}